\documentclass[reqno]{amsart}

\usepackage{amsthm}

\usepackage{%srcltx,
amssymb,amsmath,enumerate,stackrel,tabu,mathrsfs,enumitem,verbatim,stmaryrd,
leftidx,bbm, booktabs}
\usepackage{xcolor} % für Hintergrundfarben
\usepackage{xr-hyper}
\usepackage{nicefrac,stackrel}
\usepackage[all,cmtip]{xy}
\usepackage[utf8]{inputenc} % for BibTeX from arXiv
\chardef\cprime"7E % ditto
\usepackage{cases}

\usepackage[T1]{fontenc} % to get \k{e} in Kestutis (in bib file) right...

\usepackage[pagebackref,hyperindex,linktocpage=true]{hyperref}
\hypersetup{
    colorlinks,
    linkcolor={red!50!black},
    citecolor={red!50!black}, %{blue!50!black},
    urlcolor={red!50!black}, % {blue!80!black},
    filecolor={red!50!black} % black
}

\usepackage{tikz}%Für Bildchen
\usetikzlibrary{matrix,arrows,cd}

\usepackage[top=1in, left=1in, right=1in, bottom=1in, includeheadfoot]{geometry}

\usepackage{theoremref}
\definecolor{labelkey}{rgb}{1,0,0}

\makeatletter
\@addtoreset{equation}{section}
\makeatother

\numberwithin{equation}{section}

\newcommand{\timo}[1]{\framebox[1.1\width]{Timo: #1}}

\theoremstyle{definition}
\newtheorem{Defi}{Definition}[section] \newcommand{\defi}{\begin{Defi}} \newcommand{\xdefi}{\end{Defi}} %\crefname{Defi}{Def.}{Defs.}
\newtheorem{DefiLemm}[Defi]{Definition and Lemma} \newcommand{\defilemm}{\begin{DefiLemm}} \newcommand{\xdefilemm}{\end{DefiLemm}} %\crefname{DefiLemm}{Definition~and~Lem.}{Defs.~and~Lems.}
\newtheorem{Bsp}[Defi]{Example} \newcommand{\exam}{\begin{Bsp}} \newcommand{\xexam}{\end{Bsp}} %\crefname{Bsp}{Exam.}{Exam.}
\newtheorem{Syno}[Defi]{Synopsis} \newcommand{\syno}{\begin{Syno}} \newcommand{\xsyno}{\end{Syno}} %\crefname{Syno}{Syno.}{Syno.}
\newtheorem{Bem}[Defi]{Remark} \newcommand{\rema}{\begin{Bem}} \newcommand{\xrema}{\end{Bem}} %\crefname{Bem}{Rem.}{Rems.}
\newtheorem{Notation}[Defi]{Notation} \newcommand{\nota}{\begin{Notation}} \newcommand{\xnota}{\end{Notation}} %\crefname{Bem}{Rem.}{Rems.}

\theoremstyle{plain}
\newtheorem{Theo}[Defi]{Theorem} \newcommand{\theo}{\begin{Theo}} \newcommand{\xtheo}{\end{Theo}} %\crefname{Theo}{Thm.}{Thms.}
\newtheorem{Satz}[Defi]{Proposition} \newcommand{\prop}{\begin{Satz}} \newcommand{\xprop}{\end{Satz}} %\crefname{Satz}{Prop.}{Props.}
\newtheorem{Lemm}[Defi]{Lemma} \newcommand{\lemm}{\begin{Lemm}} \newcommand{\xlemm}{\end{Lemm}} %\crefname{Lemm}{Lem.}{Lems.}
\newtheorem{Coro}[Defi]{Corollary} \newcommand{\coro}{\begin{Coro}} \newcommand{\xcoro}{\end{Coro}}
\newtheorem{Ques}[Defi]{Question} \newcommand{\ques}{\begin{Ques}} \newcommand{\xques}{\end{Ques}}
\newtheorem{Conj}[Defi]{Conjecture} \newcommand{\conj}{\begin{Conj}} \newcommand{\xconj}{\end{Conj}}

\newcommand{\refsect}[1]{Section \ref{sect--#1}}

\newcommand{\eqn}{\begin{equation}} \newcommand{\xeqn}{\end{equation}}
\newcommand{\eqnarr}{\begin{eqnarray*}} \newcommand{\xeqnarr}{\end{eqnarray*}}
\newcommand{\eqnarra}{\begin{eqnarray}} \newcommand{\xeqnarra}{\end{eqnarray}}

\newcommand{\pf}{\begin{proof}} \newcommand{\xpf}{\end{proof}}

\numberwithin{equation}{section}

\newcommand{\nc}{\newcommand}%Faulheit
\nc{\StP}[1]{\cite[\href{http://stacks.math.columbia.edu/tag/#1}{Tag #1}]{StacksProject}}\nc{\StPd}[3]{\cite[Tags~\href{http://stacks.math.columbia.edu/tag/#1}{#1},
\href{http://stacks.math.columbia.edu/tag/#2}{#2},
\href{http://stacks.math.columbia.edu/tag/#3}{#3}]{StacksProject}} % Stacks-Project

\nc{\on}{\operatorname}
\nc{\aff}{{\on{aff}}}
\nc{\modi}{{\on{mod}}} % modified
\nc{\even}{{\on{even}}}
\nc{\odd}{{\on{odd}}}
\nc{\naive}{{\on{naive}}}
\nc{\hofib}{\on{hofib}}
\nc{\Bun}{\on{Bun}}
\nc{\ad}{{\on{ad}}}
\nc{\lft}{{\on{lft}}}

\nc{\Weil}{{\on{Weil}}} 
\nc{\FWeil}{{\on{FWeil}}} 

\nc{\str}{\on{-}}
\nc{\perf}{{\on{perf}}}
\nc{\Rel}{{\on{Pos}}}
\nc{\lan}{\langle}
\nc{\ran}{\rangle}
\nc{\tw}[1]{\langle #1 \rangle} % twist for Weil sheaves

\nc{\bbA}{{\mathbb A}} %% affine space
\nc{\bbB}{{\mathbb B}}
\nc{\bbC}{{\mathbb C}}
\nc{\bbD}{{\mathbb D}}
\nc{\bbE}{{\mathbb E}}
\nc{\bbF}{{\mathbb F}}
\nc{\bbG}{{\mathbb G}}
\nc{\bbH}{{\mathbb H}}
\nc{\bbI}{{\mathbb I}}
\nc{\bbJ}{{\mathbb J}}
\nc{\bbK}{{\mathbb K}}
\nc{\bbL}{{\mathbb L}}
\nc{\bbM}{{\mathbb M}}
\nc{\bbN}{{\N}} %% natural numbers
\nc{\bbO}{{\mathbb O}}
\nc{\bbP}{{\mathbb P}} % projective space
\nc{\bbQ}{{\mathbb Q}} %rational numbers
\nc{\bbR}{{\mathbb R}}
\nc{\bbS}{{\mathbb S}}
\nc{\bbT}{{\mathbb T}}
\nc{\bbU}{{\mathbb U}}
\nc{\bbV}{{\mathbb V}}
\nc{\bbW}{{\mathbb W}}
\nc{\bbX}{{\mathbb X}}
\nc{\bbY}{{\mathbb Y}}
\nc{\bbZ}{{\mathbb Z}}

\nc{\calA}{{\mathcal A}}
\nc{\calB}{{\mathcal B}}
\nc{\calC}{{\mathcal C}}
\nc{\calD}{{\mathcal D}}
\nc{\calE}{{\mathcal E}}
\nc{\calF}{{\mathcal F}}
\nc{\calG}{{\mathcal G}}
\nc{\calH}{{\mathcal H}}
\nc{\calI}{{\mathcal I}}
\nc{\calJ}{{\mathcal J}}
\nc{\calK}{{\mathcal K}}
\nc{\calL}{{\mathcal L}}
\nc{\calM}{{\mathcal M}}
\nc{\calN}{{\mathcal N}}
\nc{\calO}{{\mathcal O}}
\nc{\calP}{{\mathcal P}}
\nc{\calQ}{{\mathcal Q}}
\nc{\calR}{{\mathcal R}}
\nc{\calS}{{\mathcal S}}
\nc{\calT}{{\mathcal T}}
\nc{\calU}{{\mathcal U}}
\nc{\calV}{{\mathcal V}}
\nc{\calW}{{\mathcal W}}
\nc{\calX}{{\mathcal X}}
\nc{\calY}{{\mathcal Y}}
\nc{\calZ}{{\mathcal Z}}

\nc{\bone}{{\mathbbm{1}}}

\nc{\Sht}{{\on{Sht}}}
\nc{\Frob}{{\on{Frob}}}
\nc{\Hecke}{{\on{Hecke}}}
\nc{\inv}{{\on{inv}}}
\nc{\Conv}{{\on{Conv}}}
\nc{\triv}{{\on{triv}}}
\nc{\Isom}{{\on{Isom}}}

\nc{\scrB}{{\mathscr{B}}}
\nc{\scrA}{{\mathscr{A}}}
\nc{\bbf}{{\mathbf{f}}}
\nc{\bba}{{\mathbf{a}}}
\nc{\rig}{{\mathrm rig}}

\nc{\al}{\alpha}
\nc{\be}{\beta}
\nc{\ga}{\gamma}
\nc{\la}{\lambda}
\nc{\qcqs}{{\on{qcqs}}}
\nc{\Bmu}{{\mbox{$\raisebox{-0.59ex}{$l$}\hspace{-0.18em}\mu\hspace{-0.88em}\raisebox{-0.98ex}{\scalebox{2}{$\color{white}.$}}\hspace{-0.416em}\raisebox{+0.88ex}{$\color{white}.$}\hspace{0.46em}$}{}}}

\nc{\pot}[1]{ [\hspace{-0,5mm}[ {#1} ]\hspace{-0,5mm}] }
\nc{\rpot}[1]{ (\hspace{-0,7mm}( {#1} )\hspace{-0,7mm}) }

\nc{\defined}{\hspace{0.1cm}\stackrel{\text{\tiny \rm def}}{=}\hspace{0.1cm}}
\nc{\co}{\colon}
\nc{\specto}{{\leadsto}}

\newcommand{\category}[1]{\mathrm{#1}}
\newcommand{\Sp}{\category{Sp}} % spectra
\def\Gm{\mathbf {G}_\mathrm m} % multiplicative group
\def\SL{\mathrm {SL}} % additive group
\def\PGL{\mathrm {PGL}} % projective linear gp.
\def\GL{\mathrm {GL}} % additive group
\font\tencyr=wncyr10
\font\sevencyr=wncyr7
\font\fivecyr=wncyr5
\def\To#1#2{\mathop{\count0=#1 \loop\ifnum\count0>0 \smash-\mkern-7mu \advance\count0 -1 \repeat \mathord\rightarrow}\limits^{#2}} % long arrow with a label
\def\Hom{\mathop{\rm Hom}\nolimits} % internal hom
\def\Gr{\mathop{\rm Gr}\nolimits} % grassmannian
\def\Sht{\mathop{\rm Sht}\nolimits} % sthukas
\definecolor{hellgrau}{RGB}{200,200,200} % hellgrau
\definecolor{dunkelgrau}{RGB}{160,160,160} % grau
\definecolor{hellblau}{RGB}{194, 215, 249} %
\definecolor{dunkelblau}{RGB}{68, 128, 226} %

\def\Z{{\mathbb{Z}}} % integers and the free simplicial abelian group
\def\N{{\bf N}} % natural numbers
\def\Q{{\mathbb Q}} % rationals
\def\C{{\bf C}} % complex numbers
\def\Gm{\mathbf {G}_\mathrm m} % multiplicative group
\def\der{{\rm der}} % derived group
\def\scon{{\rm sc}} % simply connected group

\def\Spec{\mathop{\rm Spec}} % spectrum
\newcommand{\diag}{\mathrm{diag}} % operad of diagrams
\def\sbuildrel#1\over#2{\mathrel{\smash{\mathop{\kern0pt #2}\limits^{#1}}}}

\let\Gets\longleftarrow
\let\x\times

\newcommand{\lr}{\longrightarrow}

\newcommand{\SO}{\mathrm{SO}}
\newcommand{\PSO}{\mathrm{PSO}}
\nc{\supp}{\mathrm{supp}}
\nc{\cw}[1]{\omega_{#1}^\vee}
\nc{\Mlamu}{M_{\la, \mu}}

\begin{document}

\title[]{Normality of Schubert varieties \\ in affine Grassmannians}
\author[P.\,Bieker, T.\,Richarz]{Patrick Bieker and Timo Richarz}

\address{ Fakultät für Mathematik, Universität Bielefeld, Postfach 100 131, 33501 Bielefeld, Germany}
\email{pbieker@math.uni-bielefeld.de}

\address{Fachbereich Mathematik, TU Darmstadt, Schlossgartenstrasse 7, 64289 Darmstadt, Germany}
\email{richarz@mathematik.uni-darmstadt.de}

\thanks{The first named author P.B.~acknowledges support by the Deutsche Forschungsgemeinschaft (DFG, German Research Foundation) TRR 358 \textit{Integral Structures in Geometry and Representation Theory}, project number 49139240.
The second named author T.R.~is funded by the European Research Council (ERC) under Horizon Europe (grant agreement nº 101040935), by the Deutsche Forschungsgemeinschaft (DFG, German Research Foundation) TRR 326 \textit{Geometry and Arithmetic of Uniformized Structures}, project number 444845124 and the LOEWE professorship in Algebra, project number LOEWE/4b//519/05/01.002(0004)/87.}

\maketitle

\begin{abstract} 
We classify all normal Schubert varieties in the affine Grassmannian of a semisimple group over an arbitrary field with special attention to small positive characteristic.
The proof is elementary and relies on tangent space calculations for quasi-minuscule Schubert varieties, a refined Levi lemma in positive characteristic and the classification of minimal degenerations.  
\end{abstract}

%\tableofcontents

\thispagestyle{empty}

\section{Introduction}
Schubert varieties in classical flag varieties are known to be normal by work of Ramanan--Ramanathan \cite{RamananRamanathan:Schubert}, Seshadri \cite{Seshadri:Schubert}, Andersen \cite{Andersen:Schubert} and Mehta--Srinivas \cite{MehtaSrinivas:Schubert}. 
The normality of Schubert varieties in Kac--Moody affine flag varieties is due to Kumar \cite{Kumar:Demazure}, Mathieu \cite{Mathieu:CharacterFormel} and Littelmann \cite{Littelmann:Modules}. 
For affine flag varieties, interpreted as moduli spaces of torsors over an algebraic disc together with a trivialization over the punctured one, Faltings \cite{Faltings:Loops} proves normality of Schubert varieties for all split, simply connected groups. 
More generally, Pappas--Rapoport \cite{PappasRapoport:Twisted} and Lourenço \cite{Lorenco:Normality} show normality whenever the characteristic of the ground field is either zero or big enough in comparison to the group. 

For a semisimple group $G$ over an algebraically closed field $k$, a special case of \cite[Theorem 0.3]{PappasRapoport:Twisted} shows normality of Schubert varieties in the affine Grassmannian $\Gr_G$ whenever ${\rm char}(k)\nmid \#\pi_1(G)$ where $\pi_1(G)$ is the algebraic fundamental group of $G$.
For almost simple groups $G$ the number $\#\pi_1(G)$ divides the connection index of the Dynkin type of $G$ from \cite[Tables]{Bourbaki:Lie456}, and agrees with it whenever $G$ is simple.
The normality results along with other techniques are the basis for finer geometric properties such as having rational singularities, and are used in geometric representation theory or in the reduction to positive characteristic of Shimura varieties. 

When ${\rm char}(k)\mid \#\pi_1(G)$ the work of Haines, Louren\c{c}o and the second named author \cite{HainesLourencoRicharz:Normality} shows that almost all Schubert varieties in $\Gr_G$ are non-normal.
For example, if ${\rm char}(k)\mid n$, then only finitely many Schubert varieties in the affine Grassmannian for $G=\PGL_{n,k}$ are normal. 
The results of \cite{HainesLourencoRicharz:Normality} allow to classify all normal Schubert varieties for groups of small rank such as $G=\PGL_{2,k}$ for ${\rm char}(k)=2$, or $G=\PGL_{3,k}$ for ${\rm char}(k)=3$. 
In the present manuscript we classify normal Schubert varieties in $\Gr_G$ for all semisimple groups $G$. 

\subsection{Results} 
Let $G$ be an almost simple reductive group over an algebraically closed field $k$, i.e. a semisimple group with connected Dynkin diagram.
%Following \cite{HainesLourencoRicharz:Normality} in this paper a reductive group is called almost simple if its Dynkin diagram is connected.
Fix a maximal torus $T$ contained in a Borel subgroup $B$ in $G$.
Let $X_*(T)^+$ be the monoid of dominant cocharacters equipped with the Bruhat partial order.
For each $\mu\in X_*(T)^+$ the Schubert variety
\[
\Gr_{G,\leq\mu} \defined \overline{L^+G\cdot \mu(\varpi)\cdot L^+G/L^+G},
\]
is an orbit closure in the affine Grassmannian $\Gr_G$, see \refsect{Schubert-schemes}. 
By \cite{Faltings:Loops, PappasRapoport:Twisted} the $k$-variety $\Gr_{G,\leq\mu}$ is normal for all $\mu\in X_*(T)^+$ whenever ${\rm char}(k)\nmid \#\pi_1(G)$.
On the other hand if ${\rm char}(k)\mid \#\pi_1(G)$, then $\Gr_{G,\leq\mu}$ is normal for only finitely many elements $\mu\in X_*(T)^+$ by \cite{HainesLourencoRicharz:Normality}.
Note that $X_*(T)^+$ is an infinite set whenever $G$ is nontrivial. 

To state the classification consider the inclusion $X_*(T) \subset P_G^\vee$ into the coweight lattice for $G$.
Let $n\in\Z_{\geq 0}$ be the rank of $G$, and denote by $\cw1,\ldots,\cw n$ the $\Z$-basis of $P_G^\vee$ of fundamental coweights in the notation of \cite[Tables]{Bourbaki:Lie456}.
The inclusion restricts to the partially ordered monoids $X_*(T)^+ \subset (P_G^\vee)^+$ and $\cw1,\ldots,\cw n$ form a monoid basis of $(P_G^\vee)^+$.
Together with the normality result when ${\rm char}(k)\nmid \#\pi_1(G)$ the following theorem classifies all normal Schubert varieties in $\Gr_G$:

\theo[\refsect{classification}]
\thlabel{classification-normal-schubert-varieties}
Assume ${\rm char}(k)\mid \#\pi_1(G)$.
Let $\mu\in X_*(T)\backslash \{0\}$ be dominant. 
Then, the $k$-variety $\Gr_{G, \leq \mu}$ is normal if and only if the pair $(G,\mu)$ belongs to the following list:
\begin{itemize}
\item any $G$, and $\mu$ minuscule \textup{(}for a list, see \cite[Section 5, Table 1]{HePappasRapoport:semistable}\textup{)}; 
\item $G$ of type $A_n$ for $n\geq 2$ \textup{(}so ${\rm char}(k)\mid n+1$\textup{)}, and $\mu\leq d\cw 1$ or $\mu\leq d\cw n$ for some $d\in\{2,\ldots,n\}$;
\item $G$ of type $D_n$ for $n \geq 4$ \textup{(}so ${\rm char}(k) = 2$\textup{)}, and 
\begin{itemize}
	\item 
	$G \simeq \mathrm{PSO}_{2n,k}$ for $n$ odd and $\mu\in\{\cw1+\cw{n-1}, \cw1+\cw{n}\}$, or
	\item
	$G \not \simeq \mathrm{SO}_{2n,k}, \mathrm{PSO}_{2n,k}$ for $n\equiv 2\mod 4$ and $\mu = \cw1+\cw{n-1}$;	
\end{itemize}
\item $G$ adjoint of type $E_6$ \textup{(}so ${\rm char}(k)=3$\textup{)}, and $\mu\in\{2 \cw 1, \cw  3, \cw 5, 2\cw 6\}$;
\item $G$ adjoint of type $E_7$ \textup{(}so ${\rm char}(k)=2$\textup{)}, and $\mu=\cw 2$.
\end{itemize}
\xtheo

In particular, the classification is independent of the characteristic of $k$ with ${\rm char}(k)\mid \#\pi_1(G)$ and outside type $D$ also independent of the isomorphism type of $G$. 
For a list of all simple groups $G$ with ${\rm char}(k)\mid \#\pi_1(G)$ the reader is referred to \cite[Section 6]{HainesLourencoRicharz:Normality} or \refsect{classification}.
In type $D_n$ for $n$ even there is a fourth isomorphism class of semisimple groups apart from ${\rm Spin}_{2n}$, ${\rm SO}_{2n}$ and ${\rm PSO}_{2n}$. This is the group appearing in the second part of the assertion in type $D_n$ in the theorem. We identify its root datum with the one in \cite{Bourbaki:Lie456} in such a way that $\cw n$ is a cocharacter of this group.

The proof of \thref{classification-normal-schubert-varieties} is elementary in that no cohomological and representation-theoretic methods are used. 
In fact, most cases are handled by the following proposition which extends \cite[Corollary 6.2]{HainesLourencoRicharz:Normality} relying on tangent space calculations. 
Recall that the quasi-minuscule element $\mu_{\rm qm}\in X_*(T)^+$ is the unique element $\mu \in X_*(T)^+$ that is minimal with the property $\mu >0$.

\prop[\refsect{almost-minuscule}]\thlabel{prop:intro}
Assume ${\rm char}(k)\mid \#\pi_1(G)$.
Let $\la,\mu\in X_*(T)^+$ with $\la$ minuscule and $\la+\mu_{\rm qm}\leq \mu$.
Then, $\Gr_{G,\leq \mu}$ is non-normal. 
\xprop

In other words each connected component of $\Gr_G$ contains small, non-normal Schubert varieties and all bigger Schubert varieties are non-normal as well. 
The remaining cases are classified using the transversal slice and a refinement in positive characteristic of the Levi lemma from Malkin--Ostrik--Vybornov \cite{MalkinOstrikVybornov:MinimalDegenerations}, see Sections \ref{sect--transversal-slice} and \ref{sect--levi-lemma}.
The refinement in \thref{levi-lemma-derived} is necessarily more complicated as normality of Schubert varieties is sensitive to the isogeny type of $G$.  
The methods lead in \thref{mindeg-nonnormal} to the classification of non-normal minimal degeneration singularities, which is of independent interest to us. 
Even though many Schubert varieties are non-normal when ${\rm char}(k)\mid \#\pi_1(G)$ it turns out that most minimal degeneration singularities are still normal, see \thref{mindeg-nonnormal}.
So, somewhat surprisingly the combination of \thref{prop:intro} with the study of minimal degeneration singularities is sufficient to prove \thref{classification-normal-schubert-varieties} in \refsect{classification}. 

Regarding the normal locus of Schubert varieties in Dynkin type $A$ we have the following result where we denote by $\Delta$ the set of $B$-simple roots in $G$ with respect to $T$:

\coro[{\refsect{classification--typeA}}]\thlabel{coro:intro}
Let $G$ be of type $A_n$ with ${\rm char}(k)\mid \#\pi_1(G)$.
For $\mu\in X_*(T)^+$ let $J_\mu$ be the subset of all $\la \in X_*(T)^+$ such that $\mu-\la=\sum_{\al\in \Delta}n_\al\al^\vee$ for some $n_\al\in \Z_{\geq 0}$ and with at least one $n_\al=0$.
Then, the open subvariety $\cup_{\la\in J_\mu}\Gr_{G,\la}$ inside $\Gr_{G,\leq \mu}$ is normal. 
\xcoro

The exact determination of the normal locus seems challenging. 
In fact, the open locus in \thref{coro:intro} agrees with the normal locus in all examples known to us.    

\medskip
\noindent\textbf{Acknowledgements.} 
It is our pleasure to thank Tom Haines and João Lourenço for comments on an early draft.
We also thank the anonymous referee for their thoughtful suggestions, which improved the quality of the manuscript.

%%%%%%%Preliminaries on affine Grassmannians%%%%%%%%
\section{Schubert varieties}
Let $G$ be a connected reductive group over an algebraically closed field $k$.
For a $k$-algebra $R$ we denote by $R\pot{\varpi}$, respectively $R\rpot{\varpi}$ the ring of formal power series, respectively Laurent series in the formal variable $\varpi$.
The {\it loop group} $LG$ (respectively, {\it positive loop group} or {\it jet group} $L^+G$) is the group-valued functor on the category of $k$-algebras $R$ defined by $LG(R)=G(R\rpot{\varpi})$ (respectively, $L^+G(R)=G(R\pot{\varpi})$). 
Then, $L^+G\subset LG$ is a subgroup functor. 
The affine Grassmannian is the \'etale quotient
\eqn
\Gr_{G}\defined LG/L^+G,
\xeqn
which is representable by an ind-projective ind-scheme over $k$. 
The affine Grassmannian is equipped with a left action by $L^+G$.
One has $L^+G=\lim_{i\geq 0} G_i$ where $G_i(R)=G(R[\varpi]/(\varpi^{i+1}))$ is the group of $i$-jets.
Each $G_i$ is a smooth, affine $k$-group scheme with connected fibers.

\subsection{Definition}\label{sect--Schubert-schemes}
Fix a maximal torus $T$ contained in a Borel subgroup $B$ in $G$.
We denote by $\lan\str,\str\ran$ the natural paring between the characters $X^*(T)$ and cocharacters $X_*(T)$.
The sum of the $B$-positive roots is denoted $2\rho\in X^*(T)$.
Let $(X_*(T)^+,\leq)$ be the partially ordered monoid of $B$-dominant cocharacters equipped with the Bruhat partial order: 
for $\la,\mu \in X_*(T)^+$ one has $\la\leq\mu$ if and only if $\mu-\la$ is a sum of $B$-positive coroots with $\Z_{\geq 0}$-coefficients. 

The map $\mu\mapsto \mu(\varpi)$ identifies $X_*(T)^+$ with a subset of $\Gr_G(k)$. 
For each $\mu\in X_*(T)^+$, the Schubert variety $\Gr_{G,\leq \mu}$ (respectively, {$L^+G$-orbit} $\Gr_{G,\mu}$) is defined as the scheme-theoretic image (respectively, sheaf-theoretic image) of the orbit map $L^+G\to \Gr_{G}$, $g\mapsto g\cdot \mu(\varpi)$.
Then, $\Gr_{G, \leq\mu}$ is representable by a projective $k$-variety of relative dimension $\lan2\rho,\mu\ran$.
Each orbit $\Gr_{G,\mu}$ is representable by a smooth $k$-variety.
The inclusion $\Gr_{G,\mu}\subset \Gr_{G,\leq \mu}$ is open.
On the underlying sets one has $\Gr_{G,\leq \mu}=\sqcup_{\la \in X_*(T)^+: \la \leq \mu} \Gr_{G,\la}$. 
The left-$L^+G$-action on each Schubert variety factors through the group of $i$-jets $L^+G\to G_{i}$ for some $i\geq 0$.

\subsection{Normality}\label{sect--Schubert-schemes-normality}
We denote by $\pi_1(G)$ the algebraic fundamental group of $G$ which is isomorphic to the quotient of $X_*(T)$ by the coroot lattice. 
In particular, $\pi_1(G)$ is a finite group whenever $G$ is semisimple. 
The following theorem shows that most Schubert varieties are normal:

\theo[{\cite[Theorem 0.3]{PappasRapoport:Twisted}}]
\thlabel{theo:normal}
Assume either ${\rm char}(k)=0$, or ${\rm char}(k)>0$ and ${\rm char}(k)\nmid \#\pi_1(G_\der)$ where $G_\der$ is the derived group of $G$.
Then, the Schubert variety $\Gr_{G,\leq \mu}$ is normal for all $\mu\in X_*(T)^+$.
\xtheo

In addition, normal Schubert varieties are known to be Cohen--Macaulay with only rational singularities. 
The proof of \thref{theo:normal} is by reduction to the simply connected cover $G_\scon\to G_\der$, which is an étale isogeny if and only if ${\rm char}(k)\nmid \#\pi_1(G_\der)$.
The remaining cases where ${\rm char}(k)\mid \#\pi_1(G_\der)$ reduce to almost simple
%, semisimple 
groups. 
In contrast to \thref{theo:normal} almost all Schubert varieties are non-normal in these cases:

\theo[{\cite[Theorem 1.1]{HainesLourencoRicharz:Normality}}]
\thlabel{theo:not-normal}
Assume $G$ is almost simple
%, semisimple 
and ${\rm char}(k)\mid \#\pi_1(G)$.
Then, the Schubert variety $\Gr_{G,\leq \mu}$ is normal for only finitely many $\mu\in X_*(T)^+$.
\xtheo

The non-normal Schubert varieties are known to be geometrically unibranch and regular in codimension $1$, but neither S2 nor seminormal. 
The methods in \cite{HainesLourencoRicharz:Normality} allow to give an explicit bound on the number of normal Schubert varieties in $\Gr_G$ based on the following result:

\prop[{\cite[Corollary 6.2]{HainesLourencoRicharz:Normality}}]
\thlabel{prop:quasi-minuscule}
Assume $G$ is almost simple
%, semisimple 
and $p\mid \#\pi_1(G)$.
Then, the Schubert variety $\Gr_{G,\leq \mu_{\rm qm}}$ is non-normal where $\mu_{\rm qm}\in X_*(T)^+$ is quasi-minuscule. 
\xprop

%%%%%%%%%%%Minuscule plus quasi-minuscule is not normal%%%%%%%%%%
\section{Minuscule plus quasi-minuscule}\label{sect--almost-minuscule}
In this section we prove \thref{prop:intro}. 
We continue to fix $T\subset B\subset G$ as in \refsect{Schubert-schemes}.
In what follows the elements in $X_*(T)^+$ that are minimal with respect to the Bruhat partial order $\leq$ are called {\it minuscule}.
In particular, $\mu=0$ is minuscule by this definition.  

Let $I\subset L^+G$ be the {\it Iwahori subgroup} defined as the preimage of $B\subset G$ under the map $L^+G\to G$ induced by $\varpi\mapsto 0$. 
To handle the $I$-orbit stratification on $\Gr_G$ we refer to Haines--Rapoport \cite[Appendix]{PappasRapoport:Twisted} or \cite{Richarz:Iwahori} for the following properties. 
Let $N_G(T)$ be the normalizer of $T$ in $G$.
The {\it Iwahori-Weyl group} (or, {\it extended affine Weyl group}) $W:=N_G(T)(k\rpot{\varpi})/T(k\pot{\varpi})$ is equipped with the structure of a quasi-Coxeter group.
The inclusion $N_G(T)(k\rpot{\varpi})\subset G(k\rpot{\varpi})$ induces a bijection
\eqn
\label{eqn:iwahori-decomposition}
W\overset{\cong}{\lr} I(k)\backslash LG(k) /I(k). 
\xeqn
One has a semidirect product decomposition $W=T(k\rpot{\varpi})/T(k\pot{\varpi})\rtimes W_0$ where $W_0$ is the Weyl group of $(G,T)$. 
Furthermore, the inclusion $X_*(T)^+\subset T(k\rpot{\varpi})$, $\mu\mapsto \mu(\varpi)$ induces a bijection 
\eqn
\label{eqn:dominant-bijection}
X_*(T)^+\overset{\cong}{\lr} W_0\backslash W/W_0.
\xeqn

\lemm
\thlabel{lem:length-zero-minuscule}
Let $\Omega\subset W$ be the subset of length zero elements. 
Then, the map $\omega\mapsto W_0\omega W_0$ induces under the inverse of \eqref{eqn:dominant-bijection} a bijection between $\Omega$ and the minuscule elements in $X_*(T)^+$.
\xlemm
\pf
By \cite[Corollary 1.8]{Richarz:Schubert} the bijection \eqref{eqn:dominant-bijection} is compatible with the Bruhat partial orders. 
This immediately implies the lemma.
\xpf

For a minuscule element $\la\in X_*(T)^+$ we denote by $\omega_\la\in \Omega$ the associated element under \thref{lem:length-zero-minuscule}.
For any element $w\in W$ we fix a lift $\dot{w}\in N_G(T)(k\rpot{\varpi})$ and denote by 
\eqn
\label{eqn:translation}
\Gr_{G}\overset{\dot{w}\cdot}{\longrightarrow} \Gr_{G}
\xeqn
the automorphism induced by left multiplication with $\dot{w}$ on the loop group $LG$. 
For a closed subscheme $X\subset \Gr_{G}$ let $\dot{w}(X)$ be the image of $X$ under \eqref{eqn:translation}.
Clearly, if $X$ is $I$-invariant, then $\dot{w}(X)$ is independent of the choice of lift $\dot{w}\in N_G(T)(k\rpot{\varpi})$ in $w$.

\lemm
\thlabel{lem:translation-by-minuscule}
Let $\la,\mu\in X_*(T)^+$ with $\la$ minuscule.
Then, $\dot{\omega}_\la^{-1}(\Gr_{G,\leq \mu+\la})$ is an $I$-invariant, closed subvariety of $\Gr_{G}$ containing $\Gr_{G,\leq\mu}$.
\xlemm
\pf
The variety $\dot{\omega}_\la^{-1}(\Gr_{G,\leq \mu+\la})$ is $I$-invariant because $\dot{\omega}I\dot{\omega}^{-1}=I$ for all $\omega\in\Omega$.
It contains $\Gr_{G,\leq\mu}$ if and only if one has
\eqn
\label{eqn:translation-containment}
\dot{\omega}_\la(\Gr_{G,\leq \mu})\subset \Gr_{G,\leq \mu+\la}.
\xeqn
As $\omega_\la$ and $\la(\varpi)$ define the same double coset by $W_0$ there exist elements $w, w'\in W_0$ such that $w\omega_\la w'=\la(\varpi)$ inside $W$.
The Schubert varieties $\Gr_{G,\leq\mu}$ and $\Gr_{G,\leq \mu+\la}$ are $L^+G$-invariant by definition. 
So \eqref{eqn:translation-containment} holds if and only if $\la(\varpi)(\Gr_{G,\leq \mu})\subset \Gr_{G,\leq \mu+\la}$.
Clearly, this containment holds true because $\la$ and $\mu$ are dominant. 
%because the convolution morphism $\Gr_{G_k,\leq \la}\tilde\times \Gr_{G_k,\leq \mu}\to \Gr_{G_k}$ factors through $\Gr_{G_k,\leq \mu+\la}$.
\xpf

For the proof we also need the following result from \cite[Corollary 2.2]{HainesLourencoRicharz:Normality}:

\lemm
\thlabel{lemm:continuity}
Let $X\subset Y$ be $I$-invariant, closed subvarieties in $\Gr_G$.   
If $X$ is non-normal, so is $Y$.
\xlemm
\pf
Any $I$-invariant, closed subvariety in $\Gr_G$ is an $I$-orbit closure of some element in $W/W_0\subset \Gr_G(k)$.
So, the lemma is a special case of \cite[Corollary 2.2]{HainesLourencoRicharz:Normality}. 
\xpf
	
\pf[Proof of \thref{prop:intro}]
Since $G$ is almost simple
%, semisimple 
and ${\rm char}(k)\mid\#\pi_1(G)$ the Schubert variety $\Gr_{G,\leq\mu_{\rm qm}}$ is non-normal by \thref{prop:quasi-minuscule}.
It is contained in $\dot{\omega}_\la^{-1}(\Gr_{G,\leq \la+\mu_{\rm qm}})$ by \thref{lem:translation-by-minuscule}. 
The latter is therefore non-normal as well by \thref{lemm:continuity} and so is $\Gr_{G,\leq \la+\mu_{\rm qm}}$.
Since any $L^+G$-invariant, closed subvariety of $\Gr_G$ is also $I$-invariant, the Schubert variety $\Gr_{G,\leq \mu}$ is non-normal for all $\mu\in X_*(T)^+$ with $\la+\mu_{\rm qm}\leq \mu$.
This proves the proposition. 
\xpf

%%%%%%%%%%%%%%Transversal slice%%%%%%%%%%%%%
\section{Transversal slices}\label{sect--transversal-slice}
The technique is useful for studying the singularities arising when passing from a stratum to another one inside a Schubert variety. 
It is used in \cite{MalkinOstrikVybornov:MinimalDegenerations, JuteauWilliamson:Kumar, KamnitzerWebsterWeekesYacobi:slices} over $\bbC$, \cite[Section 2.3]{Zhu:Introduction} over $k$ and \cite[Section 4]{Lourenco:Grassmannian} over $\Z$. 
Due to the existence of non-normal Schubert varieties some care is needed. 
For example, the analogue of \cite[Lemma 2.6]{MalkinOstrikVybornov:MinimalDegenerations} may fail in positive characteristic.

Fix $T\subset B\subset G$ as in \refsect{Schubert-schemes}.
The {\it negative loop group} $L^-G$ is the group functor on the category of $k$-algebras $R$ defined by $L^-G(R)=G(R[\varpi^{-1}])$. 
The {\it strictly negative loop group} is the group valued functor
\eqn
L^{--}G\defined \ker(L^-G\overset{\varpi^{-1}\mapsto 0}{\lr} G).
\xeqn
By \cite[Corollary 3]{Faltings:Loops} and \cite[Corollaire 4.2.11]{Lourenco:Grassmannian} (see also \cite[Theorem 2.3.1]{dHL:Frobenius}) the composition $L^{--}G\subset LG \to \Gr_{G}$ is representable, affine and an open immersion.
There is an open cover
\eqn\label{standard-open-cover}
	\Gr_{G}=\bigcup_{\mu\in X_*(T)}\mu(\varpi)\cdot L^{--}G,
\xeqn
see \cite[Definition 5 ff.]{Faltings:Loops}, \cite[Lemma 3.1]{GoertzHaines:JordanHoelder} or \cite[Equation (4.2.24)]{Lourenco:Grassmannian}.
%In particular, the quotient map $LG\to \Gr_{G}$ admits sections Zariski locally. 

\defi\thlabel{transversal-slice-defi}
For $\la,\mu\in X_*(T)^+$ the \textit{transversal slice of $\Gr_{G,\leq \mu}$ at $\la(\varpi)$} is the étale subsheaf of $\Gr_{G}$ defined as
\eqn\label{transversal-slice-defi:eq}
	\Gr_{G,\leq \mu}^{\la}\defined L^{--}G\cdot\la(\varpi)\cap \Gr_{G,\leq \mu}.
\xeqn
\xdefi

The following lemma is analogous to \cite[Section 2.5]{MalkinOstrikVybornov:MinimalDegenerations} over the complex numbers.
It is essentially contained in \cite[Lemma 6.7]{Mueller:MinimalDegenerations}, \cite[Section 2.3]{Zhu:Introduction}, \cite[Section 4]{Lourenco:Grassmannian} or \cite[Lemma 6.1.3]{dHL:Frobenius}.

\prop\thlabel{transversal-slice-lemm}
For $\la,\mu\in X_*(T)^+$ the sheaf $\Gr_{G,\leq \mu}^{\la}$ is representable by a locally closed subscheme of $\Gr_{G}$. 
It is non-empty if and only if $\la\leq\mu$.
In this case, the following properties hold:
\begin{enumerate}
	\item\label{transversal-slice-lemm:smoothly-equivalent}
		The scheme $\Gr_{G,\leq \mu}^{\la}$ is smoothly equivalent to the open subvariety $\cup_{\nu\in X_*(T)^+:\la\leq \nu\leq \mu}\Gr_{G,\nu}$ in $\Gr_{G,\leq\mu}$.
		More precisely, for each $i\geq 0$ such that the $L^+G$-action on $\Gr_{G,\leq \mu}$ factors through the group of $i$-jets $G_{i}$ the maps in the diagram
		\eqn\label{transversal-slice-lemm:smoothly-equivalent:eq1}
			\Gr_{G,\leq \mu}^{\la}\overset{\; {\rm pr}}{\Gets} G_{i}\x \Gr_{G,\leq \mu}^{\la} \overset{{\rm act}\;}{\lr} \cup_{\nu\in X_*(T)^+:\la\leq \nu\leq \mu}\Gr_{G,\nu}
		\xeqn
		are smooth and surjective.
		Here ${\rm pr}\colon (g,x)\mapsto x$ and ${\rm act}\colon (g,x)\mapsto g\cdot x$ are the projection and the action, respectively. 
	\item\label{transversal-slice-lemm:basics}
		The scheme $\Gr_{G,\leq \mu}^{\la}$ is affine, integral, geometrically unibranch, regular in codimension $1$ and of dimension $\lan2\rho,\mu-\la\ran$.
		It is normal if and only if $\cup_{\nu\in X_*(T)^+:\la\leq \nu\leq \mu}\Gr_{G,\nu}$ is so.
		This happens, for example, if ${\rm char}(k)\nmid \#\pi_1(G_\der)$ by \thref{theo:normal}.
\end{enumerate}
\xprop
\pf
As étale sheaves $L^{--}G\cdot\la(\varpi)=L^{--}G/(L^{--}G\cap \la(\varpi)L^+G\la(\varpi)^{-1})$, which is representable by a closed sub-ind-scheme of $L^{--}G$, see \cite[Proposition 2.3.3~(1)]{Zhu:Introduction} and \cite[Equation (4.2.26)]{Lourenco:Grassmannian}.
%NTS: In the notation of loc. cit. we have $L^-\calG_a=L^{--}G\rtimes B^\opp$. As $\la$ is dominant, $\la(\varpi)^{-1}B^\opp\la(\varpi)\subset L^+G$, so $L^-\calG_a\cdot \la(\varpi)=L^{--}G\cdot \la(\varpi)$. In particular, $\Gr_{G,\leq \mu}^\la$ agrees with the Richardson semi-cell for $\la$ and $\mu$ in loc. cit.
Using that $L^{--}G$ is ind-affine the transversal slice \eqref{transversal-slice-defi:eq} is representable by an affine, locally closed subscheme of $\Gr_{G}$.

The argument for the non-emptiness criterion is as follows.
If $\la\leq \mu$, then clearly $\la(\varpi)\in \Gr_{G,\leq \mu}^{\la}(k)$.
Conversely, if $\Gr_{G,\leq \mu}^{\la}$ is non-empty, then $L^{--}G\cdot\la(\varpi)\cap \Gr_{G,\nu}$ is non-empty for some $\nu\in X_*(T)^+$ with $\nu \leq \mu$. 
By \cite[Equation (4.2.31)]{Lourenco:Grassmannian} we get $\la\leq \nu\leq \mu$.
For the rest of the proof we assume $\la\leq \mu$.

\eqref{transversal-slice-lemm:smoothly-equivalent}:
The maps in \eqref{transversal-slice-lemm:smoothly-equivalent:eq1} are surjective using the non-emptiness of $L^{--}G\cdot\la(\varpi)\cap \Gr_{G,\nu}$ for all $\nu \in X_*(T)^+$ with $\la\leq \nu \leq \mu$ to see the surjectivity of the map labelled ``act''.
Since $G_i$ is $k$-smooth, the smoothness of ``pr'' follows by base change. 
The map ``act'' arises by base change along the inclusion $\cup_{\nu\in X_*(T)^+:\la\leq \nu\leq \mu}\Gr_{G,\nu}\subset \Gr_{G}$ from the bottom map in the following commutative diagram
\eqn\label{transversal-slice-lemm:smoothly-equivalent:eq2}
\begin{tikzpicture}[baseline=(current  bounding  box.center)]
\matrix(a)[matrix of math nodes, 
row sep=1.5em, column sep=2em, 
text height=1.5ex, text depth=0.45ex] 
{L^+G \x L^{--}G & LG & LG \\ 
L^+G\x L^{--}G\cdot \la(\varpi) & & \Gr_{G} \\}; 
\path[->](a-1-1) edge node[above] {{\rm act}}  (a-1-2);
\path[->](a-1-2) edge node[above] {$\cdot \la(\varpi)$} node[below] {$\cong$} (a-1-3);
\path[->](a-2-1) edge node[above] {{\rm act}}  (a-2-3);
\path[->](a-1-1) edge node[left] {} (a-2-1);
\path[->](a-1-3) edge node[left] {} (a-2-3);
\end{tikzpicture}
\xeqn
where the vertical maps are the quotient maps.
As all schemes in \eqref{transversal-slice-lemm:smoothly-equivalent:eq1} are of finite presentation over $k$, it suffices to show that the bottom arrow satisfies the lifting property for formal smoothness and for first order nil thickenings of Artinian local rings, see \StP{0DYF}. 
On such rings we can lift points along the vertical quotient maps by their Zariski local triviality, see \eqref{standard-open-cover}.
The composition of the top arrows in \eqref{transversal-slice-lemm:smoothly-equivalent:eq2} is an open immersion by the discussion above \eqref{standard-open-cover}.
So it is formally smooth.
The right vertical arrow in \eqref{transversal-slice-lemm:smoothly-equivalent:eq2} is formally smooth because it is an $L^+G$-torsor. 
So the bottom arrow has the same lifting property for nil thickenings of Artinian local rings.
%NTS: Let $S\subset S'$ be a first order nil thickening of spectra of artinian local rings. Let $S\to L^+G\x L^{--}G\cdot \la(\varpi)$ be a map whose composition with ``act'' lifts to $S'\to \Gr_G$. Choose a preimage of $S\to L^+G\x L^{--}G\cdot \la(\varpi)$ under the left hand quotient map. As all the other maps are formally smooth we can lift this to a map from $S'$, providing the desired lift.

\eqref{transversal-slice-lemm:basics}:
We have already seen that $\Gr_{G,\leq\mu}^\la$ is an affine, locally closed subscheme of $\Gr_{G,\leq \mu}$. 
As the properties ``normal'', ``reduced'' and ``regular in codimension $1$'' are local in the smooth topology \StP{034D}, it follows from \eqref{transversal-slice-lemm:smoothly-equivalent} that $\Gr_{G,\leq\mu}^\la$ is normal if and only if $\cup_{\nu\in X_*(T)^+:\la\leq \nu\leq \mu}\Gr_{G,\nu}$ is so, and from \cite[Theorem~1.1]{HainesLourencoRicharz:Normality} that $\Gr_{G,\leq\mu}^\la$ is reduced and regular in codimension $1$.
Next, we claim that $\Gr_{G,\leq\mu}^\la$ is connected of dimension $\lan2\rho,\mu-\la\ran$. 
For this we consider the $\mathbb G_{m,k}$-action on $\Gr_{G,\leq \mu}$ induced by loop rotation $\mathbb G_{m,k}\x LG\to LG$, $(x,g(\varpi))\mapsto g(x\cdot \varpi)$.
Then, on topological spaces $\Gr_{G,\leq \mu}^\la$ identifies with the fiber of $\Gr_{G,\leq \mu}^-\to \Gr_{G,\leq \mu}^0$ above $\la(\varpi)$ where $\Gr_{G,\leq \mu}^-$ is the repeller and $\Gr_{G,\leq \mu}^0$ the fixed points, and thus is connected by \cite[Corollary~1.12]{Richarz:Spaces}.
For the dimension we note that ``act'' is a surjective, smooth map between connected $k$-schemes of finite type.
So we can compute the dimension of the source as the dimension of the target plus the dimension of any fiber.
Since $\Gr_{G,\leq\la}^\la=\{\la(\varpi)\}$ the fiber ${\rm act}^{-1}(\la(\varpi))\simeq {\rm Stab}_{G_i}(\la(\varpi))$ is the stabilizer of $G_i$ at $\la(\varpi)$. 
Using $G_i/{\rm Stab}_{G_i}(\la(\varpi))\cong \Gr_{G,\la}$ one easily gets the desired dimension formula. 
It remains to show that $\Gr_{G,\leq\mu}^\la$ is irreducible and geometrically unibranch.

Now, if $G_\der$ is simply connected, then $\Gr_{G,\leq \mu}^\la$ is normal by the above discussion and hence integral by connectedness.
To treat the general case we lift $\Gr_{G,\leq \mu}^\la$ along the central extension $\tilde G\to G$ from \cite[Section 2.4]{Lourenco:Grassmannian} as follows.
Let $T_\scon$ and $B_\scon$ be the preimage of $T\cap G_\der$ respectively $B\cap G_\der$ under the simply connected cover $G_\scon\to G_\der$.
The action of $T$ by conjugation lifts to $G_\scon$, and there is the central extension
\[
1\to T_\scon \to \tilde G:=G_\scon\rtimes T\to G\to 1
\]
with $\tilde G_\der =G_\scon$.
The group $\tilde T=T_\scon\x T$ is a split maximal torus in $\tilde G$ and $\tilde B= B_\scon\rtimes T$ a Borel subgroup. 
Clearly, $\la\leq \mu$ in $X_*(T)^+$ lifts to $\tilde \la\leq \tilde\mu$ in $X_*(\tilde T)^+$.
The induced map on Schubert varieties $\Gr_{\tilde G,\leq \tilde\mu}\to \Gr_{G,\leq\mu}$ is a finite, birational, universal homeomorphism by \cite[Proposition 3.5]{HainesRicharz:Normality} which induces isomorphisms on all residue fields by \cite[Proof of Lemma~4.4]{FHLR:singularities}. 
In particular, the map agrees with the seminormalization \StP{0H3G}. 
Thus, $\Gr_{\tilde G,\leq \tilde\mu}^{\tilde\la}\to \Gr_{G,\leq \mu}^{\la}$ is a finite, birational, homeomorphism inducing isomorphisms on all residue fields, so agrees with the seminormalization as $\Gr_{\tilde G,\leq \tilde\mu}^{\tilde\la}$ is integral and normal because $\tilde G_\der$ is simply connected. 
This implies that $\Gr_{G,\leq \mu}^{\la}$ is integral and geometrically unibranch.
\xpf

%\coro\thlabel{normality-slice-criterion}
%Let $\calP$ be a property of schemes that is local in the smooth topology.
%Then $\Gr_{G_k,\leq \mu}^\la$ has property $\calP$ if and only if $\cup_{\nu\in X_*(T)^+:\la\leq \nu\leq \mu}\Gr_{G_k,\nu}$ has property $\calP$.
%\xcoro
%This applies, for example, to the property $\calP$ of being ``regular'' or ``normal'', see \StP{034D}.

%%%%%%Levi Lemma%%%%%%%%%%%%%
\section{A Levi lemma}\label{sect--levi-lemma}
In this section we prove a refinement of the Levi lemma from \cite[Section 3.3]{MalkinOstrikVybornov:MinimalDegenerations} adapted to our purposes.
The main difference is that due to the phenomenon of non-normal Schubert varieties we need to work with derived groups as opposed to adjoint groups, which complicates the situation. 

Fix $T\subset B\subset G$ as in \refsect{Schubert-schemes}.
Let $M\subset G$ be a Levi subgroup containing $T$.  
The closed immersion on affine Grassmannians $\Gr_{M}\to \Gr_{G}$ induces a locally closed immersion
\eqn\label{immersion-Levi-slice}
	\Gr_{M,\leq \mu}^{\la} \lr \Gr_{G,\leq \mu}^{\la}
\xeqn
for any $\la,\mu\in X_*(T)^+$ with $\la\leq \mu$.
Note that every $B$-dominant cocharacter in $X_*(T)$ is $B\cap M$-dominant, so the transversal slice for the Levi subgroup is well-defined, see \thref{transversal-slice-defi}.

\lemm\thlabel{levi-lemma}
If $\mu-\la$ lies in the coroot lattice of $M$, then the map \eqref{immersion-Levi-slice} is an isomorphism. 
\xlemm
\pf
The inclusion $M\subset G$ induces a closed immersion $L^{--}M\to L^{--}G$, hence a closed immersion on orbits $L^{--}M\cdot\la(\varpi)\to L^{--}G\cdot\la(\varpi)$ and so on transversal slices, compare with the proof of \thref{transversal-slice-lemm}.
By assumption on $\mu-\la$ we have $\lan2\rho_{M},\mu-\la\ran=\lan2\rho_G,\mu-\la\ran$ where $2\rho_{M}$ and $2\rho_{G}$ denotes the sum of the roots in $B\cap M$ and $B$ respectively.
So, \eqref{immersion-Levi-slice} is a closed immersion of $k$-varieties of the same dimension and hence an isomorphism.
\xpf

Let $M_\der$ denote the derived group of $M$, and $M_\ad$ its adjoint group. 
The map $M\to M_\ad$ induces a map on transversal slices that, however, might not be an isomorphism in positive characteristic.
Therefore, we are only able to pass from $M$ to $M_\der$ under certain restrictions on the cocharacters as follows.
Let $Z$ be the maximal torus in the center of $M$.
Denote by $\Bmu:= Z\cap M_\der$ the intersection.
%The formation of $Z$ and $\Bmu$ commutes with base change, see \cite[Corollary 5.3.3]{Conrad:ReductiveGroupSchemes}. 

\lemm\thlabel{central-isog-adjoint-lemm}
The map $Z\x M_\der \to M$, $(z,m)\mapsto z\cdot m$ induces an isomorphism
\eqn\label{central-isog-adjoint}
	Z\x^{\Bmu} M_\der \overset{\cong}{\lr} M,
\xeqn
and $\Bmu=Z\cap M_\der$ is a finite $k$-group scheme of multiplicative type.
In particular, if $\Bmu$ is trivial (for example, if $M_\der$ is adjoint), then $M=Z\x M_\der$.
\xlemm
\pf
The map $Z\x M_\der \to M$ is a central isogeny.
%, see \cite[Corollary 5.3.3]{Conrad:ReductiveGroupSchemes}.
As its kernel is clearly $\Bmu$, the map \eqref{central-isog-adjoint} is an isomorphism. 
Further, $\Bmu$ is a closed $k$-subgroup of the center $Z_{M_\der}$.
As $Z_{M_\der}$ is finite over $k$ of multiplicative type, so is $\Bmu$.
%, see \cite[Remark B.1.2]{Conrad:ReductiveGroupSchemes}.
\xpf

The intersection $T\cap M_\der$ is a split maximal torus of $M_\der$.
Then, $Z\x^\Bmu(T\cap M_\der)=T$ under \eqref{central-isog-adjoint}.
In particular, the induced map on cocharacter lattices
\eqn\label{cochars-central-isog-der}
	X_*(T\cap M_\der) \oplus X_*(Z) \lr X_*(T).
\xeqn 
is injective with cokernel isomorphic to $X^*(\Bmu)$, a finite abelian group of order $\textrm{dim}_{k}\Gamma(\Bmu,\calO_{\Bmu})$.
So \eqref{cochars-central-isog-der} is an isomorphism after rationalization.  
For every $\la \in X_*(T)_\Q$, let $\la_{M_\der}$ be its projection to $X_*(T\cap M_\der)_\Q$.

%comment
\begin{comment}
We usually apply the \thref{subtori-cochars} in the following form. Let $T' = G_\der \cap T$ be the corresponding maximal $k$-split torus of $G_\der$ contained in $T$.
The central isogeny $G_\der \times Z \to G$ induces a central isogeny $T' \times Z \to T$ on tori. By the previous lemma we obtain a map
\eqn\label{cochars-central-isog-der}
X_*(T') \oplus X_*(Z) \hookrightarrow X_*(T).
\xeqn 
As before, for $\la \in X_*(T)_\Q$ we denote by $\la'$ its projection to $X_*(T')_\Q$. 
\end{comment}
%end

\lemm\thlabel{cochars-derivedgroup}
A cocharacter $\la\in X_*(T)$ is in the image of \eqref{cochars-central-isog-der} if and only if $\la_{M_\der}\in X_*(T\cap M_\der)$.
This holds, for example, if $\la_{M_\der}=0$ or if $M_\der$ is adjoint.
Further, one has 
\eqn\label{cochars-derivedgroup:eq}
	X_*(T\cap M_\der) = X_*(T) \cap \bbQ\Phi^\vee(M,T),
\xeqn
where $\Phi^\vee(M,T)\subset X_*(T)$ is the set of coroots of $M$ with respect to $T$.
\xlemm
\pf
First off, for any subtorus $T'\subset T$ its cocharacter lattice $X_*(T')$ is saturated in $X_*(T)$.
So \eqref{cochars-derivedgroup:eq} follows from $X_*(T\cap M_\der)_\Q=\bbQ\Phi^\vee(M,T)$ applying the saturation property to $T'=T\cap M_\der$.
Now let $\la \in X_*(T)$.
If $\la$ lies in the image of \eqref{cochars-central-isog-der}, then clearly $\la_{M_\der}\in X_*(T\cap M_\der)$.
Conversely, if $\la_{M_\der}\in X_*(T\cap M_\der)$, then $\la - \la_{M_\der} \in X_*(Z)_\Q$ and so $\la - \la_{M_\der} \in X_*(Z)$ applying the saturation property to $T'=Z$.
\xpf

\coro[Levi lemma]
\thlabel{levi-lemma-derived}
Let $\la,\mu\in X_*(T)^+$ with $\la\leq\mu$. 
If $\mu-\la$ lies in the coroot lattice of $M$ and if $\la_{M_\der} \in X_*(T\cap M_\der)$, then the maps
\eqn\label{levi-lemma-derived:eq}
 \Gr_{M_{\der},\leq \mu_\der}^{\la_\der} \lr \Gr_{M,\leq \mu_\der}^{\la_{\der}} \xrightarrow{(\la-\la_{\der})(\varpi)\cdot} \Gr_{M,\leq \mu}^{\la} \xrightarrow{\eqref{immersion-Levi-slice}} \Gr_{G,\leq \mu}^{\la}
\xeqn
are isomorphisms where $\mu_\der:=\mu_{M_\der}$ and $\la_\der:=\la_{M_\der}$.
\xcoro
\pf
By assumption one has $\mu-\la=(\mu-\la)_\der=\mu_\der-\la_{\der}$, so $\la_{\der}\in X_*(T\cap M_\der)$ implies $\mu_{\der} \in X_*(T\cap M_\der)$ and $\la_{\der}\leq\mu_{\der}$.
Hence, all maps are well-defined.
Next, the last map in \eqref{levi-lemma-derived:eq} is an isomorphism by \thref{levi-lemma}, so is the first map by an analogous argument.
By the proof of \thref{cochars-derivedgroup}, the cocharacter $\kappa:=\la - \la_{\der} = \mu - \mu_{\der}$ is central for $M$. 
So multiplication by $\kappa(\varpi)$ induces the middle isomorphism in \eqref{levi-lemma-derived:eq}. 
\xpf

\section{Minimal degenerations}
In this section we classify the minimal degeneration singularities that are normal. 
%This relies on Stembridge's classification of 
%We use the result to classify normal Schubert varieties but the result is also of independent interest. 
We continue to fix $T\subset B\subset G$ as in \refsect{Schubert-schemes}.
%Recall the following definition from Stembridge \cite{Stembridge:Dominant}:

\defi
A pair $\la,\mu\in X_*(T)^+$ with $\la\leq \mu$ is called a \textit{minimal degeneration} if $\mu$ lies directly above $\la$, that is, $\la<\mu$ and for all $\nu\in X_*(T)^+$ with $\la\leq\nu\leq\mu$ one has either $\la=\nu$ or $\nu=\mu$.
The singularities arising in the transversal slice $\Gr_{G,\leq \mu}^\la$ are called \textit{minimal degeneration singularities}.
\xdefi

\rema
If $\la<\mu$ is a minimal degeneration in $X_*(T)^+$, then $\Gr_{G,\leq \mu}^\la$ is smoothly equivalent to $\Gr_{G,\mu}\cup\Gr_{G,\la}$, see \thref{transversal-slice-lemm}.
So, the minimal degeneration singularities capture the singularities arising at the boundary of $\Gr_{G,\mu}$ inside $\Gr_{G,\leq\mu}$.
\xrema

For an almost simple
%, semisimple 
group $G$ one has an inclusion $X_*(T)\subset P_G^\vee$ into the coweight lattice of $G$.
Let $n$ be the rank of $G$, and denote by $\cw1,\ldots,\cw n$ the $\Z$-basis of $P_G^\vee$ of fundamental coweights in the notation of \cite[Tables]{Bourbaki:Lie456}.  
The inclusion restricts to the partially ordered monoids $X_*(T)^+ \subset (P_G^\vee)^+$ and $\cw1,\ldots,\cw n$ form a monoid basis of $(P_G^\vee)^+$.
For $\lambda < \mu$ in $X_*(T)^+$ we denote by $I_{\la,\mu}$ the set of simple coroots that appear in $\mu - \lambda$ with non-zero coefficients. 
Then, $I_{\la,\mu}$ corresponds to a subset of the Dynkin diagram for $G$. 
We denote by $M_{\la,\mu} \subseteq G$ the Levi subgroup corresponding to $I_{\la,\mu}$.
As in \thref{cochars-derivedgroup} for $\la\in X_*(T)$ let $\la_{M_\der}$ be its projection to $X_*(M_\der \cap T)_\Q$ where $M=M_{\la,\mu}$.  
The following is a special case of Stembridge's classification of minimal degenerations:

\theo
Assume $G$ is almost simple
%, semisimple
and not of Dynkin type $G_2$. 
Let $\la,\mu\in X_*(T)^+$ with $\la< \mu$, and denote by $M=M_{\la,\mu}$ the associated Levi subgroup of $G$. 
Then, $\lambda < \mu$ is a minimal degeneration if and only if one of the following conditions is satisfied:
\begin{enumerate}
	\item $\mu - \lambda$ a simple coroot;
	\label{thmStemSimp}
	
	\item $\mu - \lambda$ quasi-minuscule for $M_\der$ and $\lambda_{M_\der} = 0$;
	\label{thmStemQM}
	
	\item $\mu - \lambda$ quasi-minuscule for $M_\der$, which is of type $C_n$ for some $n \geq 2$, and $\lambda_{M_\der} = \cw n$ is the $n$-th fundamental coweight of $M_\der$.	
	\label{thmStemCn}
\end{enumerate}
\thlabel{thmStembridge}
\xtheo
\pf
This is \cite[Theorem 2.8]{Stembridge:Dominant}.
Namely, the Dynkin diagram of $G_2$ does not arise as a subdiagram of a different Dynkin type in \cite[Tables]{Bourbaki:Lie456}.
Hence, the Levi subgroup $M$ has no factor of type $G_2$, so this case is excluded from the classification.
\xpf

Using \thref{thmStembridge} together with \thref{levi-lemma-derived} we obtain the following criterion for non-normal minimal degeneration singularities:

\prop\thlabel{mindeg-nonnormal}
Assume $G$ is almost simple.% and semisimple. 
Let $\lambda < \mu$ be a minimal degeneration in $X_*(T)^+$, and denote by $M:=M_{\la,\mu}$ the associated Levi subgroup. 
Then, the minimal degeneration singularities in $\Gr_{G,\leq \mu}^\lambda$ are non-normal precisely in the following cases:
\begin{enumerate}
	\item $\mu-\lambda$ a simple coroot, $M_\der \simeq \PGL_{2,k}$ and ${\rm char}(k) = 2$;
	\item \label{mindeg-nonnormal:2}
		$\mu - \lambda$ quasi-minuscule for $M_\der$, $\lambda_{M_\der} = 0$ and ${\rm char}(k) \mid \#\pi_1(M_\der)$;
	\item $\mu - \lambda$  quasi-minuscule for $M_\der$, which is of type $C_n$ for some $n \geq 2$, $\lambda_{M_\der} = \cw n$ and ${\rm char}(k) = 2$. 
\end{enumerate}
\xprop
\begin{proof}
	First off, if $G$ is of type $G_2$, then its connection index is $1$. 
	So, any minimal degeneration singularity is normal by \thref{levi-lemma} and \thref{transversal-slice-lemm}.
	Therefore, \thref{thmStembridge} applies. 
	We now go through the different cases:
	\begin{enumerate}
		\item 
		As the corresponding Levi subgroup $M$ has semisimple rank 1, the assertion in this case follows from \thref{levi-lemma} and \thref{semi-simple-rank-1}.
		\item 
		By \thref{levi-lemma-derived} and \thref{transversal-slice-lemm} the minimal degeneration singularity is smoothly equivalent to the quasi-minuscule Schubert variety for $M_\der$. 
		This variety is non-normal if and only if ${\rm char}(k) \mid \pi_1(M_\der)$ by \thref{theo:normal} and \thref{prop:quasi-minuscule}.
		\item 
		As $M_\der$ is of type $C_n$, it is either simply connected or adjoint, in which case it has connection index $2$. 
		If $M_\der$ is simply connected or ${\rm char}(k)\neq 2$, then the minimal degeneration singularity is normal by \thref{levi-lemma} and \thref{transversal-slice-lemm}.
		If $M_\der$ is adjoint and ${\rm char}(k)= 2$, \thref{levi-lemma-derived} and \thref{transversal-slice-lemm} reduce us to the case of the Schubert variety for $\mu = \cw 1+\cw n$ in the affine Grassmmannian for $M_\der \simeq \mathrm{PSP}_{2n, k}$. 
		As $\cw 1$ is quasi-minuscule and $\cw n$ is minuscule this variety is non-normal by \thref{prop:intro}.
	\end{enumerate} 
\end{proof}

\rema
	Section \ref{sect--classification} shows that the first case in \thref{mindeg-nonnormal} only appears for $G \simeq \SO_{2n+1, k}$ with $n\geq 1$, including $\SO_{3, k}\simeq \PGL_{2, k}$.
	%Namely, when the vertex corresponding to the simple coroot $\alpha = \mu - \la$ in the Dynkin diagram for $G$ has a simple edge, \thref{lemAnLevi} shows that $M_\der \cong \SL_{2,k}$. Thus, it remains to check the cases $\alpha = \alpha_n$ for $G \cong \SO_{2n+1, k}$ and $G \cong {\rm PSP}_{2n, k}$. The case $G \cong \SO_{2n+1, k}$ is treated in the proof of \thref{lem:classification-Bn} and the case $G \cong {\rm PSP}_{2n, k}$ can be shown similarly.
\xrema

\subsection{Minimal degeneration singularities in Dynkin type $A_1$}
The singularities are calculated in \cite[Lemma 5.1]{MalkinOstrikVybornov:MinimalDegenerations} and \cite[Example 2.2]{KamnitzerWebsterWeekesYacobi:slices} over $\mathbb C$.
The manuscript \cite{Mueller:MinimalDegenerations} works over $k$ of characteristic different from $2$.
%(The condition is not explicitly mentioned but used implicitly.)
Based on \cite[Appendix B]{HainesLourencoRicharz:Normality} the result is generalized in \cite{Kumar:MinimalDegenerations} to all $k$. 
To explain the result we denote by
\eqn
	A_{G,\leq \mu}^\la\defined \Gamma(\Gr_{G,\leq\mu}^\la,\calO_{\Gr_{G,\leq\mu}^\la})
\xeqn
the ring of functions of the transversal slice $\Gr_{G,\leq\mu}^\la$ which is an affine scheme.
The following is a slight generalization of \cite{Kumar:MinimalDegenerations}:

\prop\thlabel{semi-simple-rank-1}
Assume $G$ has semi-simple rank $1$.
Let $\la,\mu\in X_*(T)^+$ with $\la<\mu$ be a minimal degeneration. 
Then one of the following cases holds:
\begin{enumerate}
	\item\label{semi-simple-rank-1:sl2}
		One has $G_\der\simeq \SL_{2,k}$. 
		Then, the ring $A_{G,\leq \mu}^\la$ is isomorphic to 
		\eqn\label{semi-simple-rank-1:GL2}
			k[x,y,z]/(z^{\lan2\rho,\mu\ran}+x y),
		\xeqn
		which is a normal domain of dimension $2$.
	\item\label{semi-simple-rank-1:pgl2}
		One has $G_\der\simeq \PGL_{2,k}$.
		Then, the ring $A_{G,\leq \mu}^\la$ is isomorphic to the subring of \eqref{semi-simple-rank-1:GL2} generated by the elements $x, y, xz, yz, z^{2}, 2z$.
		In particular, if ${\rm char}(k)\neq2$ the rings $A_{G,\leq \mu}^\la$ and \eqref{semi-simple-rank-1:GL2} are isomorphic whereas the domain $A_{G,\leq \mu}^\la$ is non-normal if ${\rm char}(k)=2$.
\end{enumerate}
\xprop
\pf
It is well-known that $G\simeq Z\times G'$ where $Z$ is a torus and $G'=\GL_{2,k}$ or $G'=\SL_{2,k}$ in case \eqref{semi-simple-rank-1:sl2}, and $G'=\PGL_{2,k}$ in case \eqref{semi-simple-rank-1:pgl2}.
Thus, we may reduce to the two cases $G=\GL_{2,k}$ and $G=\PGL_{2,k}$ which are treated in \cite{Kumar:MinimalDegenerations}. 
For convenience we give a different argument here. 

{\eqref{semi-simple-rank-1:sl2}: $G=\GL_{2,k}$.} 
Fix $T$ to be the subgroup of diagonal matrices and $B$ the subgroup of upper triangular ones.
Then $X_*(T)^+\subset \Z^2$ is the subset of elements $(a,b)\in\Z^2$ with $a\geq b$.
Since $\la=(\la_1,\la_2)<\mu=(\mu_1,\mu_2)$ in $X_*(T)^+$ is a minimal degeneration we have $\mu_1=\la_1+1$ and $\mu_2=\la_2+1$.
Further, $\lan2\rho,\mu\ran=\mu_1-\mu_2$ and similarly for $\la$, so $\lan2\rho,\mu-\la\ran=2$.
Then, one checks that the map 
\eqn\label{calculation-gl2:eq1}
(x,y,z)\mapsto \left(1+\varpi^{-1}\begin{pmatrix} \sum_{i=0}^{\la_1-\la_2}(-z)^{i+1}\varpi^{-i} & (-1)^{\la_2-\la_1}x \\ y\varpi^{\la_2-\la_1} & z \\ \end{pmatrix} \right)\cdot \la(\varpi) 
\xeqn
induces a closed immersion
\eqn\label{calculation-gl2:eq2}
	\Spec\left(k[x,y,z]/(z^{\mu_1-\mu_2}+x y)\right) \lr \Gr_{\GL_{2,k},\leq\mu}^\la.
\xeqn
As both schemes are integral of dimension $2$ by \thref{transversal-slice-lemm}, the map \eqref{calculation-gl2:eq2} is an isomorphism, see also \cite[Lemma B.3]{HainesLourencoRicharz:Normality} for a similar argument. 

{\eqref{semi-simple-rank-1:pgl2}: $G=\PGL_{2,k}$.}
Under the quotient map $\GL_{2,k}\to\PGL_{2,k}$ the minimal degeneration $\la < \mu$ in $\PGL_{2,k}$ lifts to a minimal degeneration $\tilde \la < \tilde \mu$ in $\GL_{2,k}$ as in \eqref{semi-simple-rank-1:sl2}.
This induces a map on Schubert varieties $\Gr_{\GL_{2,k},\leq\tilde\mu}\to \Gr_{\PGL_{2,k},\leq\mu}$, so a map on transversal slices
\eqn\label{pgl2-map}
	\Spec(A_{\leq\tilde\mu}^{\tilde\la}):=\Gr_{\GL_{2,k},\leq\tilde\mu}^{\tilde\la}\lr \Gr_{\PGL_{2,k},\leq\mu}^{\la}=: \Spec(A_{\leq \mu}^{\la}).
\xeqn
The map \eqref{pgl2-map} identifies $A_{\leq\tilde\mu}^{\tilde\la}$ with the normalization of the domain $A_{\leq \mu}^{\la}$ by the proof of \thref{transversal-slice-lemm}~\eqref{transversal-slice-lemm:basics}.
It remains to identify the subring $A_{\leq \mu}^{\la}$ under the isomorphism $A_{\leq\tilde\mu}^{\tilde\la}\cong k[x,y,z]/(z^m+xy)$ from \eqref{calculation-gl2:eq2} where $m:=\lan2\rho,\mu\ran$. 
For this consider the adjoint representation $\GL_{2,k}\to \GL(\mathfrak{gl}_{2,k})=\GL_{4,k}$ 
%using the ordered basis $\left (\begin{smallmatrix} 1 & 0 \\ 0 & 1 \\ \end{smallmatrix} \right )$, $\left (\begin{smallmatrix} 0 & 1 \\ 0 & 0 \\ \end{smallmatrix} \right )$, $\left (\begin{smallmatrix} 1 & 0 \\ 0 & -1 \\ \end{smallmatrix} \right )$, $\left (\begin{smallmatrix} 0 & 1 \\ 0 & 0 \\ \end{smallmatrix} \right )$, $\left (\begin{smallmatrix} 0 & 0 \\ 1 & 0 \\ \end{smallmatrix} \right )$ of $\mathfrak{gl}_{2,k}$.
explicitly given by the formula
\eqn\label{adjoint-rep}
	\begin{pmatrix}
		a & b \\
		c & d
	\end{pmatrix} 
	\mapsto \frac{1}{ad-bc} 
	\begin{pmatrix}
		ad-bc & 0 & 0 & 0\\
		0 &ad+bc & -ac & bd\\
		0 & -2ab & a^{2} & -b^{2} \\
		0 & 2cd & -c^{2} & d^{2} \\
	\end{pmatrix}.
\xeqn
This induces a closed immersion $\PGL_{2,k}\to \GL_{4,k}$, hence a closed immersion $\Gr_{\PGL_{2,k}}\to \Gr_{\GL_{4,k}}$ of affine Grassmannians.
In particular, the image of $\Gr_{\PGL_{2,k},\leq\mu}$ in $\Gr_{\GL_{4,k}}$ agrees with the scheme theoretic image of $\Gr_{\GL_{2,k},\leq\tilde\mu}$ under $\Gr_{\GL_{2,k}}\to \Gr_{\GL_{4,k}}$.
In order to describe $A_{\leq \mu}^{\la}$ as a subring of $A_{\leq \tilde\mu}^{\tilde\la}$ one computes the image of the matrix in \eqref{calculation-gl2:eq1} under \eqref{adjoint-rep}. 
Then $A_{\leq \mu}^{\la}$ is the subring of $\Z[x,y,z]/(z^{m}+xy)$ generated by the coefficients, which lie in $\Z[x,y,z]$, as polynomials in $\varpi^{\pm 1}$ of all entries. 
The diagonal matrix $\tilde\la(\varpi)$ in \eqref{calculation-gl2:eq1} maps under the adjoint representation to $\diag(1,1,\varpi^{n},\varpi^{-n})$ where $n:=m-2=\lan2\rho,\la\ran$. 
As multiplication by this matrix does not effect the coefficients as polynomials in $\varpi^{\pm 1}$ we only have to compute the image of
\eqn\label{calculation-pgl2:eq1}
	1+\varpi^{-1}
	\begin{pmatrix} 
		\sum_{i=0}^{n}(-z)^{i+1}\varpi^{-i} & (-1)^{-n}x \\ 
		y\varpi^{-n} & z \\ 
	\end{pmatrix}
	=
	\begin{pmatrix} 
		1+\sum_{i=0}^{n}(-\frac{z}{\varpi})^{i+1} & (-1)^{-n}x\varpi^{-1} \\ 
		y\varpi^{-n-1} & 1+z\varpi^{-1} \\
	\end{pmatrix}
\xeqn
under the adjoint representation. 
From the entries $bd$, $-ac$ and $d^2$ in \eqref{adjoint-rep} we see that $x,xz$, $y,yz$ and $2z, z^2$ respectively lie in $A_{\leq \mu}^{\la}$, so do $x^2$ and $y^2$ covering the entries $-b^2$ and $-c^2$.  
As $2A_{\leq \tilde\mu}^{\tilde\la}$ lies in the subring generated by $x, xz, y, yz, 2z, z^2$ the coefficients of the remaining entries in \eqref{adjoint-rep} lie in it as well.
\xpf

%%%%%%%%%Classification%%%%%%%%%%%%%%	
\section{The classification}
\label{sect--classification}
In this section we finish the proof of \thref{classification-normal-schubert-varieties}, that is, we classify all normal Schubert varieties in characteristic ${\rm char}(k)\mid \#\pi_1(G)$.
As the root systems of type $E_8$, $F_4$ and $G_2$ have connection index $1$, these cases are excluded from our list. 

Let $T \subset B \subset G$ be as in \refsect{Schubert-schemes} with $G$ almost simple
%, semisimple 
and ${\rm char}(k)\mid \#\pi_1(G)$.
We use the notation of the previous sections and \cite[Tables]{Bourbaki:Lie456}:
Let $n$ be the rank of $G$. 
We denote by $\Phi^\vee(G, T)$ the set of coroots for $G$ with respect to $T$, by $\Delta = \{ \alpha_1, \ldots, \alpha_n\}$ the set of $B$-simple roots and by $\Delta^\vee=\{\alpha_1^\vee, \ldots, \alpha_n^\vee\}$ the corresponding simple coroots.
The fundamental coweights $\cw 1, \ldots, \cw n$ form the $\Z$-basis of the coweight lattice $P_G^\vee$ dual to $\alpha_1, \ldots, \alpha_n$.
For a subset $I\subset \Delta^\vee$ we denote by $M_I$ the corresponding Levi subgroup of $G$. 
As before, for $\la,\mu\in X_*(T)^+$ with $\la \leq \mu$ let $I_{\la, \mu}\subset \Delta^\vee$ be the subset of simple coroots that appear in $\mu -\la$ with a non-zero coefficient. 
In this case we also write $\Mlamu = M_{I_{\la,\mu}}$. 
We use the notation
$$\Z^n_{\text{even}} = \big\{(a_1, \ldots, a_n) \in \Z^n \colon \textstyle{\sum_{i=1}^n} a_i \equiv 0 \text{ mod } 2 \big\} \quad \text{and} \quad \Z^n_{\Sigma = 0} = \big\{(a_1, \ldots, a_n) \in \Z^n \colon \textstyle{\sum_{i=1}^n} a_i = 0 \big\}.$$

%%%%%%%%%%%%%strategy%%%%%%%%%%%%%%
\subsection{Strategy of proof}
Let $\mu_{\rm qm} \in X_*(T)^+$ be the quasi-minuscule element.
We freely use the following results during the proof: 
\begin{enumerate}
	\item \label{eqn:principle1}
		The Schubert variety $\Gr_{G,\leq \mu}$ is non-normal whenever $\la + \mu_{\rm qm} \leq \mu$ for some minuscule $\la \in X_*(T)^+$ , see \thref{prop:intro}.
	\item \label{eqn:principle2}
		One has $\Gr_{M,\leq \mu}^\la\cong\Gr_{G,\leq \mu}^\la$ for all $\la,\mu\in X_*(T)^+$ with $\la \leq \mu$ and $M=\Mlamu$, see \thref{levi-lemma}.
		Under the condition and notation in \thref{levi-lemma-derived}, one also has $\Gr_{M_\der,\leq \mu_\der}^{\la_\der}\cong \Gr_{G,\leq \mu}^\la$.
	\item \label{eqn:principle3}
		If ${\rm char}(k)\nmid \#\pi_1(M_\der)$ with $M$ as in \eqref{eqn:principle2} (for example, $M_\der$ simply connected), then $\Gr_{G,\leq \mu}^\la$ is normal and so is $\cup_{\nu\in X_*(T)^+:\la\leq\nu\leq\mu}\Gr_{G,\nu}$, see \thref{theo:normal} and \thref{transversal-slice-lemm}.
\end{enumerate}
With principles \eqref{eqn:principle1}--\eqref{eqn:principle3} in hand we use a case by case analysis relying on the classification of groups \cite[Tables]{Bourbaki:Lie456}.
In many cases $M_\der$ turns out to be simply connected.
The remaining cases are handled by \eqref{eqn:principle1} applied to $M$ and \thref{mindeg-nonnormal}.

For the convenience of the reader, we list the (quasi-)minuscule coweights in each relevant type, and its decomposition into simple coroots from \cite[Tables]{Bourbaki:Lie456}.
\begin{equation*}
\begin{tabular}{c|c|c}
	Type & quasi-minuscule coweight & (non-zero) minuscule coweights\\
	\hline 
	$A_n$, $n \geq 1$ & $ \alpha^\vee_1 + \ldots + \alpha^\vee_n$ & $\cw 1, \ldots, \cw n $ \\
	$B_n$, $n \geq 3$ &  $\cw 2 = \alpha^\vee_1 + 2(\alpha_2^\vee + \ldots + \alpha_{n-1}^\vee) + \alpha_n^\vee $ & $\cw 1$\\
	$C_n$, $n \geq 2$ & $ \cw 1 = \alpha^\vee_1 +  \alpha^{\vee}_2 + \ldots +  \alpha_{n}^{\vee} $  & $\cw n$\\ 
	$D_n$, $n \geq 4$ & $\cw 2 = \alpha_1^\vee + 2(\alpha_2^\vee + \ldots + \alpha_{n-1}^\vee) + \alpha_n^\vee $  & $\cw 1, \cw{n-1}, \cw{n}$ \\\
	$E_6$ & $\cw 2 = \alpha_1^\vee + 2 \alpha_2^\vee +  2 \alpha_3^\vee + 3\alpha_4^\vee + 2 \alpha_5^\vee + \alpha_6^\vee$ & $\cw{1}, \cw 6$ \\
	$E_7$ & $\cw 1 = 2\alpha_1^\vee + 2 \alpha_2^\vee +  3 \alpha_3^\vee + 4 \alpha_4^\vee + 3 \alpha_5^\vee + 2 \alpha_6^\vee + \alpha_7^\vee$ & $\cw{7}$
\end{tabular}
\end{equation*}

%%%%%%%%%%%%%%%type A_n%%%%%%%%%
\subsection{Type $A_n$, $n \geq 1$}
\label{sect--classification--typeA}
The case $G = \PGL_{2,k}$ is treated in \cite[Corollary 6.8]{HainesLourencoRicharz:Normality}, see also \thref{semi-simple-rank-1}. 
So, we consider $n \geq 2$. 
Let $G \simeq \SL_{n+1, k}/\Bmu_{\ell, k}$ for some $\ell \mid (n+1)$ with ${\rm char}(k) \mid \ell$. 

\lemm\thlabel{lemm:simply-connected-type-A}
Let $I$ be a proper subset of $\Delta^\vee$. 
Then, the derived subgroup of $M = M_I$ is simply connected.
In particular, when $I$ is connected $M_\der \simeq \SL_{m+1,k}$, where $m := \#I < n$. 
\thlabel{lemAnLevi}
\xlemm
\begin{proof}
	It clearly suffices to check the case $G = \PGL_{n+1,k}$.
	Using \thref{cochars-derivedgroup} we compute 
	$$X_\ast(T \cap M_\der) = X_\ast(T) \cap \Q \Phi^\vee(M,T) = \Z^{n+1}/\langle (1, \ldots, 1) \rangle \cap \Q_{\Sigma = 0}^I  = \Z_{\Sigma = 0}^I = \Z \Phi^\vee(M,T).$$
	Hence, $M_\der$ is simply connected and thus isomorphic to $\SL_{m+1,k}$ when $I$ is connected. 
\end{proof}

We use \thref{lemm:simply-connected-type-A} to prove \thref{coro:intro}:

\pf[Proof of \thref{coro:intro}]
Let $\mu \in X_*(T)^+$ and let $\la \in J_\mu$, i.e. 
$\la \in X_*(T)^+$ with $\la < \mu$ such that $\mu-\la=\sum_{\al\in \Delta}n_\al\al^\vee$ for some $n_\al\in \Z_{\geq 0}$ and with at least one $n_\al=0$. 
In particular, $\nu \in J_\mu$ for every $\la \leq \nu \leq \mu$, and hence $\cup_{\la\in J_\mu}\Gr_{G,\la} \subset \Gr_{G,\leq \mu}$ is open.

Moreover, $I_{\la, \mu}$ is a proper subset of $\Delta^\vee$. 
So by \thref{lemAnLevi} the derived subgroup of $M=\Mlamu$ is simply connected.
Hence, $\Gr_{M,\leq \mu}^\la\cong\Gr_{G,\leq \mu}^\la$ is normal.
As the subset $J_\mu$ is the union of all $\la$ as above the open subvariety $\cup_{\la\in J_\mu}\Gr_{G,\la}$ in $\Gr_{G,\leq \mu}$ is normal. 
\xpf

\coro
For $\mu\in \{d \cw 1, d \cw n\}$ the variety $\Gr_{G,\leq \mu}$ is normal for all $d\in\{1,\ldots,n\}$. 
\xcoro
\pf
It suffices to treat the case $G = \PGL_{n+1,k}$ and $\mu = d \cw 1$.
The case $\mu = d \cw n$ can then be deduced by applying the non-trivial automorphism of the Dynkin diagram. 
Let $\la = \cw d$ which is the unique minuscule coweight with $\la \leq d \cw 1$.
It is quickly verified that $I_{\la, d \cw 1} = \{\alpha_1^\vee, \ldots, \alpha_{d-1}^\vee \}$. %and $I_{\la, d \cw n} = \{\alpha_{d+1}^\vee, \ldots, \alpha_n^\vee\}$.
With the notation of \thref{coro:intro} one has $\cup_{\la\in J_\mu}\Gr_{G,\la}=\Gr_{G,\leq \mu}$ which is normal. 
\xpf

The following lemma finishes the proof of \thref{classification-normal-schubert-varieties} in Dynkin type $A_n$:

\lemm
Let $\mu \in X_\ast(T) \setminus \{0\}$ be dominant. 
Then, $\mu$ satisfies precisely one of the following two properties:
\begin{itemize}
	\item $\mu \leq d \cw 1$ or $\mu \leq d \cw n$ for some $d\in \{1,\ldots,n\}$;
	\item $\mu \geq \la + \mu_{\rm qm}$ for some $\la\in X_*(T)^+$ minuscule. 
\end{itemize}
\xlemm
\pf
Let $\la\in X_*(T)^+$ be the unique minuscule cocharacter with $\la < \mu$. 
Following the proof of the previous corollary, for $\mu$ with the first property we have $I_{\la, \mu} \neq \{ \alpha_1^\vee, \ldots, \alpha_n^\vee\}$ so the cases are mutually exclusive.
When $\la = 0$, the claim follows from the fact that the quasi-minuscule cocharacter $\mu_{\rm qm}$ is the unique minimal degeneration of $\la$. 
Let us now assume that $\la = \cw d$ for some $d\in \{1,\ldots,n\}$ and $\mu \not \geq \la + \mu_{\rm qm}$.
As $\mu_{\rm qm} = \sum_{i = 1}^n \alpha_i^\vee$, this means that $I_{\la, \mu}$ is a proper subset of $\Delta^\vee$ and hence a subset of $\{\alpha_1^\vee, \ldots, \alpha_{d-1}^\vee\}$ or $\{\alpha_{d+1}^\vee, \ldots, \alpha_{n}^\vee\}$.
In the first case we get $\mu \leq d \cw 1$ and in the second case $\mu \leq (n+1-d) \cw n$. 
\xpf

%%%%%%%%%%%%%%%%type B_n%%%%%%%%%%%%%
\subsection{Type $B_n$, $n \geq 3$}
The case $B_2 = C_2$ is treated in the next subsection. 
The root system of type $B_n$ has connection index $2$. 
So, we consider a simple group $G \simeq \SO_{2n+1, k}$ with ${\rm char}(k)=2$. 
\thref{classification-normal-schubert-varieties} in Dynkin type $B_n$ follows from the next lemma:

\lemm
\thlabel{lem:classification-Bn}
The unique minuscule $\la := \cw 1\in X_*(T)\backslash\{0\}$ admits $\la < \mu := \cw 3$ as the unique minimal degeneration, and $\Gr_{G, \leq \mu}$ is non-normal. 
In particular, every non-minuscule $\nu\in X_*(T)^+$ satisfies $\nu\geq \mu_{\rm qm}$ or $\nu\geq \mu$, so the normal Schubert varieties in $\Gr_G$ are precisely the minuscule ones.
\xlemm
\pf 
We note that $\cw 3 - \cw 1 = \alpha_2^\vee + 2(\alpha^\vee_3 + \ldots + \alpha^\vee_{n-1}) + \alpha^\vee_n$ is the quasi-minuscule coweight for the corresponding Levi subgroup $M = \Mlamu$ which is of type $B_{n-1}$.  
%so $I_{\la, \mu} = \{\alpha_2^\vee, \ldots, \alpha_{n}^\vee\}$. 
The fact that $\la < \mu$ is a minimal degeneration of $\la$ follows from \thref{thmStembridge} \eqref{thmStemQM}. 
To see that it is indeed the unique minimal degeneration of $\la$ we note that by \thref{thmStembridge} there are a priori two other possibilities for minimal degenerations. Namely,
a quasi-minuscule minimal degenerations (case \eqref{thmStemQM}) for a Levi subgroup $M_I$ for some proper subset $I$ of $I_{\la, \mu}$, which contradicts the minimality of $\la < \mu$; or a simple minimal degeneration (case \eqref{thmStemSimp}). But the minimality of $\la < \mu$ excludes all $\la + \alpha_i^\vee$ for $i \neq 1$ and $\la + \alpha_1^\vee$ is not dominant.  
%Namely, it is case \eqref{thmStemQM} of the theorem.

By \thref{cochars-derivedgroup} together with the explicit description of coroots and coweights in \cite{Bourbaki:Lie456}, we see that 
$$ X_*(T \cap M_{\der}) = X_*(T) \cap \Q\Phi^\vee(M,T) = \bigoplus_{i \in I_{\la, \mu}} \Z (\omega_M)_i^\vee = \Z^{n-1}$$
contains $\Z \Phi^\vee(M_\der, T \cap M_\der) = \Z^{n-1}_{\rm{even}}$ as an index 2 subgroup. 
In particular, $M_\der \cong \SO_{2n-1, k}$ is adjoint, and so $\Gr_{M_\der,\leq\mu_\der}^{\la_\der}\cong\Gr_{G,\leq \mu}^\la$ is non-normal by \thref{mindeg-nonnormal}~\eqref{mindeg-nonnormal:2}.
\xpf

%%%%%%%%%%%%%%%type C_n%%%%%%%%%%%%%
\subsection{Type $C_n$, $n \geq 2$}
The root system of type $C_n$ has connection index 2. 
So, we consider an adjoint group $G \simeq \mathrm{PSP}_{2n, k} = \Sp_{2n, k}/\Bmu_{2, k}$ with ${\rm char}(k)=2$. 
We include the case $n = 2$, that is, $\SO_{5, k} \simeq \mathrm{PSP}_{4, k}$.
\thref{classification-normal-schubert-varieties} in Dynkin type $C_n$ follows from the next lemma:

\lemm
Both minuscule cocharacters $\la$ of $G$ admit $\la < \la + \mu_{\rm qm}$ as the unique minimal degeneration.
In particular, the normal Schubert varieties in $\Gr_G$ are precisely the minuscule ones.
\xlemm
\pf
When $\la=0$, the assertion is obvious.
If $\la\neq 0$, then $\la + \mu_{\rm qm}$ is a minimal degeneration by \thref{thmStembridge}. 
Namely, it corresponds to case \eqref{thmStemCn} of the theorem.
The uniqueness can be shown as in the proof of \thref{lem:classification-Bn}. 
\xpf

%%%%%%%%type D_n%%%%%%%%%%%%%%
\subsection{Type $D_n$ for $n \geq 4$}
The root system $D_n$ has connection index $4$, so ${\rm char}(k) = 2$. 
Recall that the quotient of the weight lattice by the root lattice is isomorphic to $\Z/2\Z \times \Z/2\Z$ when $n$ is even and $\Z/4\Z$ when $n$ is odd.
Thus, when $n$ is odd, $G \simeq \SO_{2n, k}$ or $G \simeq \PSO_{2n, k} = \SO_{2n, k}/\Bmu_{2, k}$. 
When $n$ is even there is a third possibility.
\thref{classification-normal-schubert-varieties} in Dynkin type $D_n$ follows from the next lemma:

\lemm
For $\mu\in X_*(T)\cap \{\cw 2, \cw 2 + \cw{n-1},\cw 2 + \cw n\}$ the Schubert variety $\Gr_{G,\leq\mu}$ is non-normal.
The minscule coweight $\cw 1$ has a unique minimal degeneration $\cw 1 < \mu = (1,1,1,0, \ldots,0)$, i.e., $\mu =\cw 3 + \cw 4$ for $n = 4$ and $\mu = \cw 3$ for $n \geq 5$. 
The Schubert variety $\Gr_{G, \leq \mu}$ is non-normal.

\begin{figure}
	\begin{tikzcd}
		0 \arrow{rr}{\{1, \ldots, n\}} && \cw 2   && \cdots 
	\end{tikzcd}
	
	\begin{tikzcd}
		&&&& \\
		\cw 1 \arrow{rr}{\{2, \ldots, n\}} && \cw 3  & & \ldots
	\end{tikzcd}

	\begin{tikzcd}
		&&&&\\
		\cw{n-1} \arrow{rr}{\{1, \ldots, n-2, n\}} && \cw 1 + \cw{n} \arrow{rr}{\{2, \ldots,n-1\}}  && \cw 2 + \cw{n-1} & \ldots 
	\end{tikzcd}
	\begin{tikzcd}
		&&&&\\
		\cw{n} \arrow{rr}{\{1, \ldots, n-1\}} && \cw 1 + \cw{n-1} \arrow{rr}{\{2, \ldots,n-2, n\}}  && \cw 2 + \cw{n} & \ldots 
	\end{tikzcd}
	
	\caption{The Bruhat partial order on dominant cocharacters for $\PSO_{2n, k}$ for $n \geq 5$. 
	Each arrow is labelled by the subset of indices of simple coroots appearing in the associated minimal degeneration.}
	\label{figure-dom-cochars-Dn}
\end{figure}
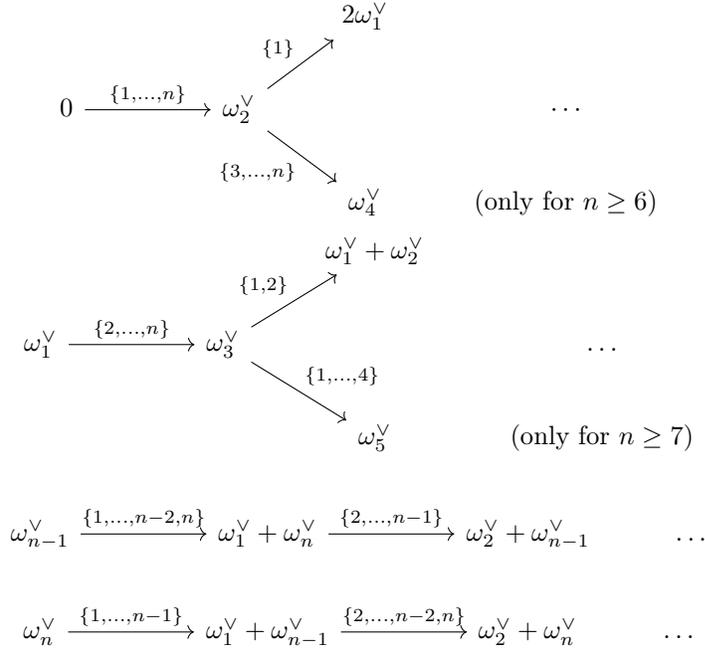

Moreover, for $\mu\in X_*(T)\cap\{\cw 1 + \cw{n-1},\cw 1 + \cw n \}$ the Schubert variety $\Gr_{G,\leq \la}$ is normal if and only if $G \simeq \PSO_{2n, k}$ for $n$ odd, or $\mu = \cw 1 + \cw{n-1}$ and  $G \not \simeq \SO_{2n, k}, \PSO_{2n, k}$ for $n\equiv 2\mod 4$.
\xlemm
\pf
First off, one can quickly verify the description of the Bruhat partial order in Figure \ref{figure-dom-cochars-Dn} using \thref{thmStembridge}.
The non-normality of the Schubert varieties for $\cw 2, \cw 2 + \cw n$ and $\cw 2 + \cw{n-1}$ follows from \thref{prop:intro}.

The coweight $\cw 1$ is a dominant cocharacter for $G \simeq \SO_{2n, k}$ and $G \simeq \PSO_{2n, k}$. 
That it admits a unique minimal degeneration $\cw 1 < \mu$ as claimed can be verified again using \thref{thmStembridge}, it is quasi-minuscule for $M_\der$ where $M=M_I$ with $I=\{\al_2^\vee, \ldots, \al_n^\vee\}$.
We have $\Z \Phi^\vee(M, T) = 0 \oplus \Z^{n-1}_{\rm{even}}$ and
$X_*(T) = \Z^n$ for $G \simeq \SO_{2n, k}$ or $X_*(T) = \Z^n + \Z(1/2, \ldots, 1/2)$ for $G \simeq \PSO_{2n, k}$. 
One has $X_\ast(T \cap M_\der) = \Q  \Phi^\vee(M, T)\cap X_\ast(T) \simeq (0 \oplus \Q^{n-1}) \cap \Z^n = \Z^{n-1}$ by \thref{cochars-derivedgroup} in both cases. 
Thus, $\pi_1(M_\der) \cong \Z/2\Z$ and the minimal degeneration singularity $\Gr^{\cw 1}_{G, \leq \mu}$ is non-normal by \thref{mindeg-nonnormal}.
Note that $M_\der \simeq \SO_{2n-2, k}$ for $n \geq 5$ and $M_\der \simeq \SL_{4, k}/\Bmu_{2, k}$ for $n = 4$. 

The coweight $\cw n$ is a cocharacter for $G \simeq \PSO_{2n, k}$ and $G \not \simeq \SO_{2n, k}, \PSO_{2n, k}$.
Its unique minimal degeneration $\cw{n} < \cw 1 + \cw{n-1}$ is quasi-minuscule for the Levi subgroup $M=M_I$ with $I = \{\al_1^\vee, \ldots, \al_{n-1}^\vee\}$. 
In this case $\Z \Phi^\vee(M, T) = \Z^n_{\Sigma = 0}$. 
For $G \simeq \PSO_{2n, k}$ using \thref{cochars-derivedgroup}, $X_*(T \cap M_\der) =  \Z^n_{\Sigma = 0} + \Z(1/2, \ldots, 1/2, -1/2, \ldots, -1/2)$, so $\#\pi_1(M_\der) =2$, when $n$ is even and 
$X_*(T \cap M_\der) =  \Z^n_{\Sigma = 0}$, so $\#\pi_1(M_\der) =1$, when $n$ is odd.
Otherwise, i.e. for $n$ even and $G \not \simeq \SO_{2n, k}, \PSO_{2n, k}$, we have $X_*(T) \cong \Z^n_{\rm even} + \Z(1/2, \ldots, 1/2)$. Hence, $X_*(T \cap M_\der) =  \Z^n_{\Sigma = 0} + \Z(1/2, \ldots, 1/2, -1/2, \ldots, -1/2)$, so $\#\pi_1(M_\der) =2$, when $n/2$ is even and 
$X_*(T \cap M_\der) =  \Z^n_{\Sigma = 0}$, so $\#\pi_1(M_\der) =1$, when $n/2$ is odd.
The claim now follows from \thref{mindeg-nonnormal}.

The coweight $\cw {n-1}$ is a cocharacter only for $\PSO_{2n, k}$, and this case can be treated analogously to the previous one.
\xpf

%%%%%%%%%%%%%%%%type e_6%%%%%%%%%%%
\subsection{Type $E_6$}
The root system of type $E_6$ has connection index 3. 
So, we consider an adjoint group $G$ of type $E_6$ with ${\rm char}(k)=3$. 
\thref{classification-normal-schubert-varieties} in Dynkin type $E_6$ follows from the next lemma:

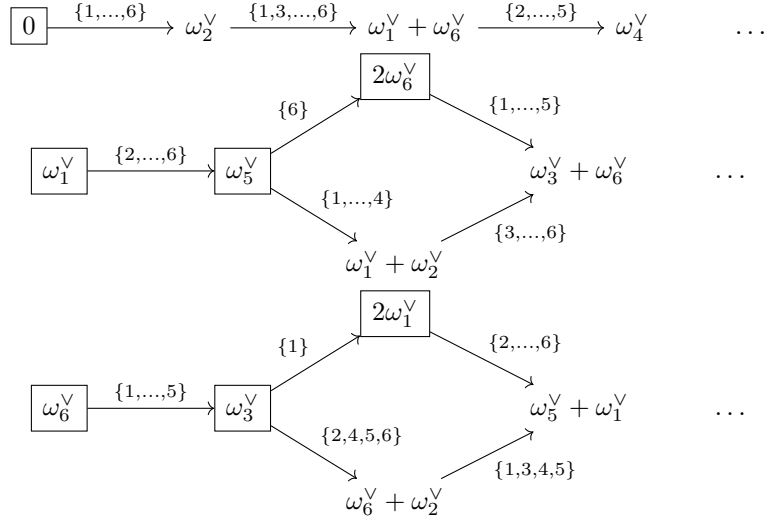
\begin{figure}
	\begin{tikzcd}
		|[draw=black]|0 \arrow{rr}{\{1, \ldots, 6\}} && \cw 2 \arrow{rr}{\{1,3,\ldots,6\}} && \cw 1 + \cw 6 \arrow{rr}{\{2,\ldots, 5\}} && \cw 4 & \ldots
	\end{tikzcd}
	
	\begin{tikzcd}
		&&& |[draw=black]|2 \cw 6 \arrow{rd}{\{1, \ldots, 5\}}  &&\\
		|[draw=black]|\cw 1 \arrow{rr}{\{2, \ldots, 6\}} && |[draw=black]|\cw 5 \arrow{ru}{\{6\}} \arrow{rd}{\{1, \ldots, 4\}} && \cw 3 + \cw 6 & \ldots \\
		&&& \cw 1 + \cw 2 \arrow{ru}[swap]{\{3, \ldots, 6\}} &&
	\end{tikzcd}
	
	\begin{tikzcd}
		&&& |[draw=black]|2 \cw 1 \arrow{rd}{\{2, \ldots, 6\}}  &&\\
		|[draw=black]|\cw 6 \arrow{rr}{\{1, \ldots, 5\}} && |[draw=black]|\cw 3 \arrow{ru}{\{1\}} \arrow{rd}{\{2,4,5,6\}} && \cw 5 + \cw 1 & \ldots \\
		&&& \cw 6 + \cw 2 \arrow{ru}[swap]{\{1,3, 4, 5\}} &&
	\end{tikzcd}
	\caption{The Bruhat partial order on dominant cocharacters for $G$ of type $E_6$. 
	The boxes indicate the normal Schubert varieties.}
	\label{figure-dom-cochars-E6}
\end{figure}

\lemm
The normal Schubert varieties in $\Gr_G$ are precisely the minuscule ones and $\Gr_{G,\leq\mu}$ for all $\mu\in\{2\cw 1, \cw 3, \cw 5, 2\cw 6\}$. 
\xlemm
\pf
As in case $D$ above one checks the Bruhat partial order in Figure \ref{figure-dom-cochars-E6} using the Stembridge classification. 
So, the Schubert varieties for all dominant cocharacters not listed in the lemma are non-normal.

Now let $\mu\in \{\cw 5,2 \cw 6\}$ and $\la = \cw 1$. 
Then, $I_{\la, \mu} = \{\al_2^\vee, \ldots, \al_6^\vee\}$. 
Hence, the Levi subgroup $\Mlamu$ is of type $D_5$. In particular, $3 \nmid \# \pi_1(M_\der)$ and so the transversal slice $\Gr^\la_{G, \leq \mu}$ is normal.
The assertion for $\mu\in\{\cw 3,2 \cw 1\}$ follows similarly.
\xpf

%%%%%%%%%%type E_7%%%%%%%%%%%
\subsection{Type $E_7$}	
The root system of type $E_7$ has connection index 2. 
So, we consider an adjoint group $G$ of type $E_7$ with ${\rm char}(k)=2$. 
\thref{classification-normal-schubert-varieties} in Dynkin type $E_7$ follows from the next lemma:
\begin{figure}
	\begin{tikzcd}%[cramped, sep=small] 
		|[draw=black]|0 \arrow{rr}{\{1,..,7\}}	& &  \cw 1 \arrow{rr}{\{2,..,7\}}	&	&  \cw 6	&		\ldots \\
	\end{tikzcd}
	
	\begin{tikzcd}%[cramped, sep=small] 
		|[draw=black]| \cw 7 \arrow{rr}{\{1,..,6\}}	&& |[draw=black]| \cw 2 \arrow{rr}{\{1,3,..,7\}}	&& \cw 1 + \cw 7 	& \dots	\\	
	\end{tikzcd}
	\caption{Same description as in Figure \ref{figure-dom-cochars-E6} but for $G$ of type $E_7$.}
	\label{figure-dom-cochars-E7}
\end{figure}
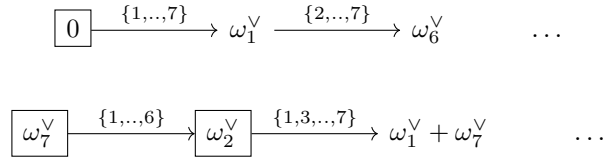

\lemm
The normal Schubert varieties in $\Gr_G$ are precisely the minuscule ones as well as $\Gr_{G,\leq \cw 2}$. 
\xlemm
\pf
We proceed as for type $E_6$. 
One checks Figure \ref{figure-dom-cochars-E7} using the Stembridge classification. 
So, the Schubert varieties for all dominant cocharacters not listed in the lemma are non-normal.

Now let $\mu = \cw 2$ and $\la = \cw 7$. 
Then, $I_{\la, \mu} = \{\al_1^\vee, \ldots, \al_6^\vee\}$. 
Hence, the Levi subgroup $\Mlamu$ is of type $E_6$. 
In particular, $2 \nmid \# \pi_1(M_\der)$ and so the transversal slice $\Gr^\la_{G, \leq \mu}$ is normal.
\xpf

\bibliographystyle{alphaurl}
\bibliography{bib}

\end{document}